\documentclass{amsart}
\usepackage{amsmath}
\usepackage{amssymb}
\usepackage{amsfonts}
\usepackage{graphicx}
\usepackage{comment}
\usepackage{epsfig}
\usepackage{epstopdf}
\usepackage{hyperref}

\newtheorem{theorem}{Theorem}[section]
\newtheorem{corollary}[theorem]{Corollary}
\newtheorem{lemma}[theorem]{Lemma}
\newtheorem{proposition}[theorem]{Proposition}
\newtheorem{definition}[theorem]{Definition}
\newtheorem{example}[theorem]{Example}
\newtheorem{remark}[theorem]{Remark}

\begin{document}
\title[Local dimensions]{Local dimensions of measures of finite type II -
Measures without full support and with non-regular probabilities. \\
\today }
\author[K.~E.~Hare, K.~G.~Hare, M.~K.~S.~Ng]{Kathryn E. Hare, Kevin G. Hare,
Michael Ka Shing Ng}
\thanks{Research of K. G. Hare was supported by NSERC Grant RGPIN-2014-03154}
\thanks{Computational support provided in part by the Canadian Foundation
for Innovation, and the Ontario Research Fund.}
\thanks{Research of K. E. Hare and M. K. S. Ng was supported by NSERC Grant
2011-44597}
\address{Dept. of Pure Mathematics, University of Waterloo, Waterloo, Ont.,
N2L 3G1, Canada}
\email{kehare@uwaterloo.ca}
\address{Dept. of Pure Mathematics, University of Waterloo, Waterloo, Ont.,
N2L 3G1, Canada}
\email{kghare@uwaterloo.ca}
\address{Dept. of Pure Mathematics, University of Waterloo, Waterloo, Ont.,
N2L 3G1, Canada}

\begin{abstract}
Consider a sequence of linear contractions $S_{j}(x)=\varrho x+d_{j}$ and
probabilities $p_{j}>0$ with $\sum p_{j}=1$. We are interested in the
self-similar measure $\mu =\sum p_{j}\mu \circ S_{j}^{-1}$, of finite type.
In this paper we study the multi-fractal analysis of such measures,
extending the theory to measures arising from non-regular probabilities and
whose support is not necessarily an interval.

Under some mild technical assumptions, we prove that there exists a subset
of supp$\mu $ of full $\mu $ and Hausdorff measure, called the truly
essential class, for which the set of (upper or lower) local dimensions is a
closed interval. Within the truly essential class we show that there exists
a point with local dimension exactly equal to the dimension of the support.
We give an example where the set of local dimensions is a two element set,
with all the elements of the truly essential class giving the same local
dimension. We give general criteria for these measures to be absolutely
continuous with respect to the associated Hausdorff measure of their support
and we show that the dimension of the support can be computed using only
information about the essential class.

To conclude, we present a detailed study of three examples. First, we show
that the set of local dimensions of the biased Bernoulli convolution with
contraction ratio the inverse of a simple Pisot number always admits an
isolated point. We give a precise description of the essential class of a
generalized Cantor set of finite type, and show that the $kth$ convolution
of the associated Cantor measure has local dimension at $x\in (0,1)$ tending
to 1 as $k$ tends to infinity. Lastly, we show that within a maximal loop
class that is not truly essential, the set of upper local dimensions need
not be an interval. This is in contrast to the case for finite type measures
with regular probabilities and full interval support.
\end{abstract}

\maketitle

\allowdisplaybreaks

\section{Introduction}

\label{sec:intro} In this paper we continue the investigations, begun in 
\cite{HHM}, on the multifractal analysis of equicontractive self-similar
measures of finite type. For self-similar measures arising from an IFS that
satisfies the open set condition the multifractal analysis is well
understood. In particular, the set of attainable local dimensions is a
closed interval whose endpoints can be computed with the Legendre transform.

For measures that do not satisfy the open set condition, the multifractal
analysis is more complicated and the set of local dimensions need not be an
interval. This phenomena was discovered first for the 3-fold convolution of
the classical Cantor measure in \cite{HL} and was further explored in \cite%
{BHM, FLW, LW, Sh}, for example. In \cite{NW}, Ngai and Wang introduced the
notion of finite type, a property stronger than the weak separation
condition (WSC), but satisfied by many interesting self-similar measures
that fail the open set condition. Examples include Bernoulli convolutions
with contraction factor the inverse of a Pisot number and self-similar
Cantor-like measures with ratio the inverse of an integer.

Building on earlier work, such as \cite{FO,Hu,LN1,Po}, Feng undertook a
study of equicontractive, self-similar measures of finite type in \cite%
{F3,F1,F2}, with his main focus being Bernoulli convolutions. Motivated by
this research, in \cite{HHM} (and \cite{Homepage}) a general theory was
developed for the local dimensions of self-similar measures of finite type
assuming the associated self-similar set was an interval and the underlying
probabilities $\{p_{j}\}_{j=0}^{m}$ generating the measure $\mu $ were
regular, meaning $p_{0}=p_{m}=\min p_{j}$. There it was shown that the set
of local dimensions at points in the `essential class' (a set of full
Lebesgue measure in the support of $\mu $ and often the interior of its
support) was a closed interval and that the set of local dimensions at
periodic points was dense in this interval. Formulas were given for the
local dimensions. These formulae are particularly simple at periodic points.

In contrast to much of the earlier work, in this paper we do not require any
assumptions on the probabilities and we relax the requirement that the
support of $\mu $ (the self-similar set) is an interval. In Section \ref%
{sec:term}, we give formulas for calculating local dimensions. These are
relatively simple for periodic points, although necessarily more complicated
than under the previous assumptions. We begin Section \ref{sec:TEP} by
introducing the `truly essential class'. We see that this set is the
relative interior of the essential class and we prove that it has full $\mu $
and Hausdorff $s$-measure, where $s$ is the Hausdorff dimension of the
self-similar set. Under a mild technical assumption, that is satisfied in
many interesting examples, we prove that the set of local dimensions at the
points in the truly essential class is a closed interval and that the set of
local dimensions at the periodic points is dense in that interval.

We prove that there is always a point at which the local dimension of $\mu $
coincides with the Hausdorff dimension of $\mathrm{supp}\mu $ and give an
example of a measure where this occurs at all the truly essential points
(but not at all points of the support). A sufficient condition is given for
a finite type measure to be absolutely continuous with respect to the
associated Hausdorff measure and an example is seen that satisfies this
condition when $s=1$, even though the self-similar set is not an interval.
We also give a formula for calculating the Hausdorff dimension of the
support from just the knowledge of the essential class.

Related results were given by Feng in \cite{F2}. Here, Feng constructed a
(typically, countably infinite) family of closed intervals, $I_{j},$ with
disjoint interiors, where $\bigcup I_{j}$ is of full measure and on each of
these closed intervals the set of attainable local dimensions of the
restricted measure $\mu _{j}:=\mu |_{I_{j}}$ was a closed interval. From his
construction one can see that $\bigcup I_{j}\bigcap K$ is the essential
class and $\bigcup \int(I_{j})\bigcap K$ is contained in our truly essential
class. The endpoints of the intervals $I_{j}$ may or may not be truly
essential. We note that the local dimension of the restricted measure $\mu
_{j}$ at an end point of $I_{j}$ is not necessarily the same as the local
dimension of $\mu $ at this point, even when it is a truly essential point.

Feng and Lau in \cite{FL} studied more general IFS that only satisfy the WSC
and showed that in this case there is also an open set $U$ such that the set
of attainable local dimensions of the restricted measure, $\mu |_{U},$ is a
closed interval. In the examples given in that paper, the set $U$ is much
smaller than the truly essential class.

In \cite{F1}, Feng had shown that the set of local dimensions of the uniform
Bernoulli convolutions with contraction factor the inverse of a simple Pisot
number (meaning the minimal polynomial is $x^{n}-\sum_{j=0}^{n-1}x^{j}$) is
always a closed interval. As one application of our work, in Section \ref%
{sec:Pisot} we prove, in contrast, that biased Bernoulli convolutions with
these contraction factors always admit an isolated point in their set of
local dimensions.

In Section \ref{sec:Cantor} we present a detailed study of the local
dimensions of finite type Cantor-like measures, extending the work done in 
\cite{BHM, HHM, Sh}. In those papers, it was shown, for example, that if $%
p_{0}<p_{j}$ for $j\neq 0,m$, then the local dimension at $0$ is isolated.
Here we give further conditions that ensure there is an isolated point. But
we also give examples where the measure has no isolated points and we give a
family of examples which have exactly two distinct local dimensions. We also
show that the local dimensions of the rescaled $k$-fold convolutions of a
Cantor-like measure converge to $1$ at points in $(0,1)$. Previously, in 
\cite{BH}, it had been shown that these local dimensions were bounded.

In Section \ref{sec:MLC}, we illustrate, by means of a detailed example, the
complications and differences that can arise when studying the local
dimensions outside of the truly essential class and in Section \ref{sec:misc}
investigate the connection between finite type and Pisot contractions.

\section{Notation and Preliminary Results}

\label{sec:term}

We begin by introducing the definition of finite type, as well as basic
notation and terminology that will be used throughout the paper.

\subsection{Finite Type}

Consider the iterated function system (IFS) consisting of the contractions $%
S_{j}:\mathbb{R\rightarrow R}$, $j=0,\dots ,m$, defined by 
\begin{equation}
S_{j}(x)=\varrho x+d_{j},
\end{equation}%
where $0<\varrho <1$, $0=d_{0}<d_{1}<d_{2}<\cdot \cdot \cdot <d_{m}$ and $%
m\geq 1$ is an integer. The unique, non-empty, compact set $K$ satisfying 
\begin{equation*}
K=\bigcup\limits_{j=0}^{m}S_{j}(K)
\end{equation*}%
is known as the \emph{associated self-similar set}. By rescaling the $d_{j}$
if needed, we can assume the convex hull of $K$ is $[0,1]$. We will not
assume that $K=[0,1]$ or even that it has non-empty interior.

It was shown in \cite[Thm. 1.2]{NW} that if $s=\dim _{H}K$ and $H^{s}$
denotes the Hausdorff $s$-measure restricted to $K$, then $0<H^{s}(K)<\infty
.$ Upon normalizing we can assume $H^{s}(K)=1$. Further, we note that $%
0<s\leq 1$. We remark that in the special case that $K=[0,1]$, then $s=1$
and $H^{s}$ is the normalized Lebesgue measure.

Suppose probabilities $p_{j}>0$, $j=0,\dots ,m$ satisfy $%
\sum_{j=0}^{m}p_{j}=1$. Throughout this paper, our interest will be in the
self-similar measure $\mu $ associated to the family of contractions $%
\{S_{j}\}$ given above, which satisfies the identity 
\begin{equation}
\mu =\sum_{j=0}^{m}p_{j}\mu \circ S_{j}^{-1}.  \label{ss}
\end{equation}%
These non-atomic, probability measures have support $K$.

We put $\mathcal{A}=\{0,\dots,m\}$. Given an $n$-tuple $\sigma
=(j_{1},\dots,j_{n})$ $\in \mathcal{A}^{n}$, we write $S_{\sigma }$ for the
composition $S_{j_{1}}\circ \cdot \cdot \cdot \circ S_{j_{n}}$ and let 
\begin{equation*}
p_{\sigma }=p_{j_{1}}\cdot \cdot \cdot p_{j_{n}}\text{.}
\end{equation*}

\begin{definition}
The iterated function system (IFS), 
\begin{equation*}
\{S_{j}(x)=\varrho x+d_{j}:j=0,\dots ,m\},
\end{equation*}%
is said to be of \textbf{finite type} if there is a finite set $F\subseteq 
\mathbb{R}$ such that for each positive integer $n$ and any two sets of
indices $\sigma =(j_{1},\dots ,j_{n})$, $\sigma ^{\prime }=(j_{1}^{\prime
},\dots ,j_{n}^{\prime })$ $\in \mathcal{A}^{n}$, either 
\begin{equation*}
\varrho ^{-n}\left\vert S_{\sigma }(0)-S_{\sigma ^{\prime }}(0)\right\vert >c%
\text{ or }\varrho ^{-n}(S_{\sigma }(0)-S_{\sigma ^{\prime }}(0))\in F,
\end{equation*}%
where $c=(1-\varrho )^{-1}(\max d_{j}-\min d_{j})$ is the diameter of $K$.

If the IFS is of finite type and $\mu $ is an associated self-similar
measure satisfying (\ref{ss}), we also say that $\mu $ is of \textbf{finite
type}.
\end{definition}

Here we have given the general definition of finite type for an
equicontractive IFS in $\mathbb{R}$. This simplifies to $c=1 $ in the case
where the convex hull of $K$ is $[0,1]$. It is worth noting here that the
definition of finite type is independent of the choice of probabilities.

Finite type is a property that is stronger than the weak separation
condition, but weaker than the open set condition \cite{Ng}. Examples
include (uniform or biased) Bernoulli convolutions with contraction factor
the reciprocal of a Pisot number, and Cantor-like measures associated with
Cantor sets with contraction factors reciprocals of integers. See Sections %
\ref{sec:Pisot} and \ref{sec:Cantor} where these are studied in detail.

\subsection{Characteristic vectors and the Essential class}

The structure of measures of finite type is explained in detail in \cite{F3,
F1, F2} and \cite{HHM}; we will give a brief overview here.

For each integer $n$, let $h_{1},\dots ,h_{s_{n}}$ be the collection of
elements of the set $\{S_{\sigma }(0),$ $S_{\sigma }(1):\sigma \in \mathcal{A%
}^{n}\}$, listed in increasing order. Put 
\begin{equation*}
\mathcal{F}_{n}=\{[h_{j},h_{j+1}]:1\leq j\leq s_{n}-1\text{ and }%
(h_{j},h_{j+1})\cap K\neq \emptyset \}\text{.}
\end{equation*}%
Elements of $\mathcal{F}_{n}$ are called \textit{net intervals of level }$n$%
. By definition, a net interval contains net subintervals of every lower
level. For each $\Delta \in \mathcal{F}_{n}$, $n\geq 1$, there is a unique
element $\widehat{\Delta }\in \mathcal{F}_{n-1}$ which contains $\Delta ,$
called the \textit{parent} (of \textit{child} $\Delta )$. We will define the 
\emph{left-most child} of parent $\widehat{\Delta }=[a,b]\in \mathcal{F}%
_{n-1}$ to be the child $\Delta =[a,b^{\prime }]\in \mathcal{F}_{n}$. We
will similarly define the \emph{right-most child}. It is worth noting that
it is possible for a child to be both the left and the right-most child. It
is further worth observing that because we are not assuming the self-similar
set is the interval $[0,1]$, it is possible for a parent to have no
left-most child or no right-most child.

Given $\Delta =[a,b]\in \mathcal{F}_{n}$, we denote the \textit{normalized
length} of $\Delta $ by 
\begin{equation*}
\ell _{n}(\Delta )=\varrho ^{-n}(b-a)\text{.}
\end{equation*}%
By the \textit{neighbour set} of $\Delta $ we mean the ordered $k$-tuple 
\begin{equation*}
V_{n}(\Delta )=(a_{1},\dots ,a_{k})
\end{equation*}%
where 
\begin{equation*}
\{a_{1},\dots ,a_{k}\}=\{\varrho ^{-n}(a-S_{\sigma }(0)):\sigma \in \mathcal{%
A}^{n}\text{, \thinspace }\Delta \subseteq S_{\sigma }[0,1]\}.
\end{equation*}%
Given $\Delta _{1},\dots ,\Delta _{m},$ (listed in order from left to right)
all the net intervals of level $n$ which have the same parent and normalized
length as $\Delta $, let $r_{n}(\Delta )$ be the integer $r$ with $\Delta
_{r}=\Delta $. The \textit{characteristic} \textit{vector of }$\Delta $ is
the triple 
\begin{equation*}
\mathcal{C}_{n}(\Delta )=(\ell _{n}(\Delta ),V_{n}(\Delta ),r_{n}(\Delta )).
\end{equation*}%
Often we suppress $r_{n}(\Delta )$ giving the \emph{reduced characteristic
vector} $(\ell _{n}(\Delta ),V_{n}(\Delta ))$.

If the measure is of finite type, there will be only finitely many distinct
characteristic vectors. We denote the set of such vectors by $\Omega ,$ 
\begin{equation*}
\Omega =\{\mathcal{C}_{n}(\Delta ):n\in \mathbb{N}\text{, }\Delta \in 
\mathcal{F}_{n}\}\text{.}
\end{equation*}

By an \textit{admissible path, }$\eta ,$\textit{\ of length }$L(\eta )=L,$
we will mean an ordered $L$-tuple, $\eta =(\gamma _{j})_{j=1}^{L},$ where $%
\gamma _{j}\in \Omega $ for all $j$ and the characteristic vector, $\gamma
_{j},$ is the parent of $\gamma _{j+1}$. Each $\Delta \in \mathcal{F}_{n}$
can be uniquely identified by an admissible path of length $n+1$, say $(%
\mathcal{C}_{0}(\Delta _{0}),\dots ,\mathcal{C}_{n}(\Delta _{n}))$, where $%
\Delta =\Delta _{n}$, $\Delta _{0}=[0,1]$, $\Delta _{j}\in \mathcal{F}_{j}$
and $\Delta _{j}=\widehat{\Delta _{j+1}}$ for all $j$. This is called the 
\textit{symbolic representation} of $\Delta ;$ we will frequently identify $%
\Delta $ with its symbolic representation.

Similarly, the \textit{symbolic representation }for $x\in K$ will mean the
sequence%
\begin{equation*}
\lbrack x]=(\mathcal{C}_{0}(\Delta _{0}),\mathcal{C}_{1}(\Delta _{1}),\dots )
\end{equation*}%
of characteristic vectors where $x\in \Delta _{n}$ for all $n$ and $\Delta
_{j}\in \mathcal{F}_{j}$ is the parent of $\Delta _{j+1}$. The notation $%
[x|N]$ will mean the admissible path consisting of the first $N$
characteristic vectors of $[x]$. We will often write $\Delta _{n}(x)$ for
the net interval in $\mathcal{F}_{n}$ containing $x\in K$; its symbolic
representation is $[x|n]$.

If $x$ is an endpoint of $\Delta _{n}(x)$ for some $n$ (and then for all
larger integers) we call $x$ a \textit{boundary point. }We remark that if $x$
is a boundary point, then there can be two different symbolic
representations for $x$, one approaching $x$ from the left i.e., by taking
right-most descendents at all levels beyond level $n,$ and the other
approaching $x$ from the right, by taking left-most descendents. If $x$ is
not a boundary point, then the symbolic representation is unique.

It is worth emphasizing that $[x|N]$ is defined as the truncation of $[x]$
as opposed to defining it as a sequence $(\mathcal{C}_{0}(\Delta _{0}),%
\mathcal{C}_{1}(\Delta _{1}),\dots, \mathcal{C}_{N}(\Delta _{N}))$ with $%
x\in \Delta _{i}$. To see this distinction, recall that it is possible for $%
\Delta _{N}=[h_{i},h_{i+1}]$ to have no right-most children. 
Let $x=h_{i+1}$ be the right-most endpoint of $\Delta
_{N} $. Then $x$ is also the left-most endpoint of the adjacent net
interval, $\Delta _{N}^{\prime }=[h_{i+1},h_{i+2}]$. As $\Delta _{N}$ has no
right-most child, we do not have a net interval of depth $N+1$ with $x\in
\Delta _{N+1}\subseteq \Delta _{N}$. As $K$ has no isolated points and $x\in
K$, for all $M\geq N$ we must have net intervals $x\in \Delta _{M}\subseteq
\Delta _{N}^{\prime }$. In such a case, the boundary point $x$ has a unique
symbolic representation.

A non-empty subset $\Omega ^{\prime }\subseteq \Omega $ is called a\textbf{\ 
}\textit{loop class} if whenever $\alpha ,\beta \in \Omega ^{\prime }$, then
there are characteristic vectors $\gamma _{j}$, $j=1,\dots ,n$, such that $%
\alpha =\gamma _{1}$, $\beta =\gamma _{n}$ and $(\gamma _{1},\dots ,\gamma
_{n})$ is an admissible path with all $\gamma _{j}\in \Omega ^{\prime }$. A
loop class $\Omega ^{\prime }\subseteq \Omega $ is called an \textit{%
essential class} if, in addition, whenever $\alpha \in \Omega ^{\prime }$
and $\beta \in \Omega $ is a child of $\alpha $, then $\beta \in \Omega
^{\prime }$. Of course, an essential class is a maximal loop class.

In \cite[Lemma 6.4]{F2}, Feng proved the important fact that there is always
precisely one essential class, which we will denote by $\Omega _{0}$. If $%
[x]=(\gamma _{0},\gamma _{1},\gamma _{2},\dots )$ with $\gamma _{j}\in
\Omega _{0}$ for all large $j$, we will say that $x$ is an \textit{essential
point} (or is\textit{\ in the essential class}) and similarly speak of a net
interval being essential. A path $(\gamma _{j})_{j=1}^{L}$ is in the
essential class if all $\gamma _{j}\in \Omega _{0}$. We similarly speak of a
point, net interval or path as being in a given loop class. The finite type
property ensures that every element in the support of $\mu $ is contained in
a maximal loop class.

We remark that the essential class is dense in the support of $\mu $. This
is because the uniqueness of the essential class ensures that every net
interval contains a net subinterval in the essential class. In Proposition %
\ref{SizeofTE} we will show that the essential class has full $\mu $ measure
and full Hausdorff $s$-measure in $K$, where $s$ is the Hausdorff dimension
of $K$.

\subsection{Transition matrices}

A very important concept in the multifractal analysis of measures of finite
type are the so-called transition matrices. These are defined as follows:
Let $\Delta =[a,b]$ be a net interval of level $n$ with parent $\widehat{%
\Delta }=[c,d]$ . Assume $V_{n}(\Delta )=(a_{1},\dots ,a_{N})$ and $V_{n-1}(%
\widehat{\Delta })=(c_{1},\dots ,c_{M})$. The \textit{primitive transition
matrix}, $T(\mathcal{C}_{n-1}(\widehat{\Delta }),$ $\mathcal{C}_{n}(\Delta
)),$ is a $M\times N$ matrix whose $jk$ entry is given by%
\begin{equation*}
T_{jk}:=\left( T(\mathcal{C}_{n-1}(\widehat{\Delta }),\mathcal{C}_{n}(\Delta
))\right) _{jk}=p_{\ell }\text{ }
\end{equation*}%
if $\ell \in \mathcal{A}$ and there exists $\sigma \in \mathcal{A}^{n-1}$
with $S_{\sigma }(0)=c-\varrho ^{n-1}c_{j}$ and $S_{\sigma \ell
}(0)=a-\varrho ^{n}a_{k}$, and $T_{jk}=0$ otherwise. We note that in \cite%
{HHM} the transition matrices are normalized so that the minimal non-zero
entry is $1$. That is, we used $p_{\ast }^{-1}T$ instead of $T$, where $%
p_{\ast }=\min p_{j}$.

We observe that each column of a primitive transition matrix has at least
one non-zero entry. The same is true for each row if supp$\mu =[0,1]$, but
not necessarily otherwise, see Example \ref{Ex:notpos}.

Given an admissible path $\eta =(\gamma _{1},\dots ,\gamma _{n})$, we write 
\begin{equation*}
T(\eta )=T(\gamma _{1},\dots ,\gamma _{n})=T(\gamma _{1},\gamma _{2})\cdot
\cdot \cdot T(\gamma _{n-1},\gamma _{n})
\end{equation*}%
and refer to such a product as a \textit{transition matrix}. We will say the
transition matrix $T(\gamma _{1},\dots ,\gamma _{n})$ is \textit{essential}
if all $\gamma _{j}$ are essential characteristic vectors.

By the norm of a matrix $T$ we mean%
\begin{equation*}
\left\Vert T\right\Vert =\sum_{jk}\left\vert T_{jk}\right\vert .
\end{equation*}%
A matrix is called \textit{positive} if all its entries are strictly
positive. An admissible path $\eta $ is called \textit{positive} if $T(\eta
) $ is a positive matrix. Here is an elementary lemma which shows the
usefulness of positivity.

\begin{lemma}
\label{compare} Assume $A,B,C$ are transition matrices and $B$ is positive.

\begin{enumerate}
\item There are constants $a,b>0$, depending on the matrices $A$ and $B$
respectively, so that $\left\Vert AC\right\Vert \geq a\left\Vert
C\right\Vert $ and $\left\Vert ABC\right\Vert \geq b\left\Vert A\right\Vert
\left\Vert C\right\Vert $. \label{it2:1}

\item If each row of $A$ has a non-zero entry, then there is a constant $\,c$%
, depending on matrix $C$, such that $\left\Vert AC\right\Vert \geq
c\left\Vert A\right\Vert $. \label{it2:2}

\item There is a constant $C_{1}=C_{1}(B)$ such that if $AB$ is a square
matrix, then 
\begin{equation*}
sp(AB)\leq \left\Vert AB\right\Vert \leq C_{1}sp(AB)\text{.}
\end{equation*}%
\label{it2:3}

\item Suppose $B$ is a square matrix. There is a constant $C_{2}=C_{2}(B)$
such that 
\begin{equation*}
sp(B^{n})\leq \left\Vert B^{n}\right\Vert \leq C_{2}sp(B^{n})\text{ for all }%
n\text{.}
\end{equation*}
\label{it2:4}
\end{enumerate}
\end{lemma}

\begin{proof}
Parts \eqref{it2:1} and \eqref{it2:2} follow by simply writing the
expressions for $\left\Vert AC\right\Vert $ and $\left\Vert ABC\right\Vert $
in terms of the entries of $A,B,C$, and noting that a transition matrix has
non-negative entries and each column has a non-zero entry.

Parts \eqref{it2:3} and \eqref{it2:4} follow as in \cite[Lemma 3.15]{HHM}.
\end{proof}

\subsection{Basic facts about of local dimensions of measures of finite type}

\begin{definition}
Given a probability measure $\mu $, by the \textbf{upper local dimension} of 
$\mu $ at $x\in $ supp$\mu $ we mean the number 
\begin{equation*}
\overline{\dim}_{loc}\mu (x)=\limsup_{r\rightarrow 0^{+}}\frac{\log \mu
([x-r,x+r])}{\log r}.
\end{equation*}%
Replacing the $\limsup $ by $\liminf $ gives the \textbf{lower local
dimension}, denoted $\underline{\dim}_{loc}\mu (x)$. If the limit exists, we
call the number the \textbf{local dimension} of $\mu $ at $x$ and denote
this by $\dim_{loc}\mu (x)$.
\end{definition}

It is easy to see that 
\begin{equation}
\dim_{loc}\mu (x)=\lim_{n\rightarrow \infty }\frac{\log \mu ([x-\varrho
^{n},x+\varrho ^{n}])}{n\log \varrho }\text{ for }x\in \text{supp}\mu ,
\label{locdim}
\end{equation}%
and similarly for the upper and lower local dimensions.

\textbf{Notation: }Throughout the paper, when we write $F_{n}\sim G_{n}$ we
mean there are positive constants $c_{1},c_{2}$ such that 
\begin{equation*}
c_{1}F_{n}\leq G_{n}\leq c_{2}F_{n}\text{ for all }n.
\end{equation*}

To calculate local dimensions, it will be helpful to know $\mu (\Delta )$
for net intervals $\Delta $.

\begin{proposition}
\label{HHM}Let $\Delta _{n}=[a,b]\in \mathcal{F}_{n},$ with $V_{n}(\Delta
_{n})=(a_{1},\dots ,a_{N})$. Then 
\begin{equation*}
\mu (\Delta _{n})=\sum_{i=1}^{N}\mu \lbrack a_{i},a_{i}+\ell _{n}(\Delta
_{n})]\sum_{\substack{ \sigma \in \mathcal{A}^{n}  \\ \varrho
^{-n}(a-S_{\sigma }(0))=a_{i}}}p_{\sigma }.
\end{equation*}%
Furthermore, if $[\Delta _{n}]=(\gamma _{0},\gamma _{1},\dots,\gamma _{n})$
and 
\begin{equation*}
P_{n}(\Delta _{n})=\sum_{i=1}^{N}\sum_{\sigma \in \mathcal{A}^{n}:\varrho
^{-n}(a-S_{\sigma }(0))=a_{i}}p_{\sigma }\text{ ,}
\end{equation*}%
then%
\begin{eqnarray*}
\mu (\Delta _{n}) &\sim &P_{n}(\Delta _{n})=\left\Vert T(\gamma _{0},\gamma
_{1},\dots,\gamma _{n})\right\Vert \\
&=&p_{\ast }^{n}\left\Vert T^{\ast }(\gamma _{0},\gamma _{1},\dots,\gamma
_{n})\right\Vert
\end{eqnarray*}%
where $p_{\ast }=\min p_{j}$.
\end{proposition}

\begin{proof}
This follows in a similar fashion to Lemma 3.2, Corollary 3.4 and the
discussion prior to Corollary 3.10 of \cite{HHM}, noting that $\mu \lbrack
a_{i},a_{i}+\ell _{n}(\Delta )]\geq \mu (S_{\sigma }^{-1}(a,b))>0$.
\end{proof}

The analogue of Proposition \ref{HHM} was very useful in \cite{HHM} as it
was the key idea in proving the following formula.

\begin{corollary}
\label{corlocdim} \cite[Cor. 3.10]{HHM} Suppose $\mu $ is a self-similar
measure satisfying identity (\ref{ss}), that has support $[0,1],$ is of
finite type and has probabilities satisfying $p_{0}=p_{m}=\min p_{j}$. If $%
x\in $supp$\mu $, then 
\begin{eqnarray}
\dim _{loc}\mu (x) &=&\frac{\log p_{0}}{\log \varrho }+\lim_{n\rightarrow
\infty }\frac{\log \left\Vert T^{\ast }([x|n])\right\Vert }{n\log \varrho }
\label{F2} \\
&=&\lim_{n\rightarrow \infty }\frac{\log \left\Vert T([x|n])\right\Vert }{%
n\log \varrho }
\end{eqnarray}%
and similarly for the upper and lower local dimensions.
\end{corollary}

This corollary need not be true, however, if the assumptions of supp$\mu
=[0,1]$ and regular probabilities, i.e., $p_{0}=p_{m}=\min p_{j},$ are not
all satisfied. Instead, we proceed as follows.

\textbf{Terminology: }Assume $\{h_{j}\}=\{S_{\sigma }(0),S_{\sigma
}(1):\sigma \in \mathcal{A}^{n}\}$ with $h_{j}<h_{j+1}$ and suppose $\Delta
_{n}=[h_{i},h_{i+1}]$ is a net interval of level $n$. Let $\Delta _{n}^{-}$
be the empty set if $(h_{i-1},h_{i})\cap K$ is empty and otherwise, let $%
\Delta _{n}^{-}$ $=$ $[h_{i-1},h_{i}]$ $.$ Similarly, define $\Delta
_{n}^{+} $ to be the net interval immediately to the right of $\Delta _{n}$
(or the empty set), with the understanding that if $\Delta _{n}$ is the left
or right-most net interval in $\mathcal{F}_{n}$, then $\Delta _{n}^{-}$
(respectively, $\Delta _{n}^{+})$ is the empty set. We refer to $\Delta
_{n}^{-}(x),\Delta _{n}(x),\Delta _{n}^{+}(x)$ as \textit{adjacent net
intervals} (even if some are the empty set).

If $x$ belongs to the interior of $\Delta _{n}(x)$, we put%
\begin{equation}
M_{n}(x)=\mu (\Delta _{n}(x))+\mu (\Delta _{n}^{+}(x))+\mu (\Delta
_{n}^{-}(x)).  \label{Mn}
\end{equation}%
If $x$ is a boundary point of $\Delta _{n}(x)=[h_{i},h_{i+1}]$, we put 
\begin{equation}
M_{n}(x)=\mu (\Delta _{n}(x))+\mu (\Delta _{n}^{\prime }(x)),  \label{Mnbdy}
\end{equation}%
where $\Delta _{n}^{\prime }(x)=\Delta _{n}^{-}(x)$ if $x=h_{i}$ and $\Delta
_{n}^{\prime }(x)=\Delta _{n}^{+}(x)$ if $x=h_{i+1}$. We will refer to $%
\Delta _{n}^{\prime }(x)$ as the \textit{other net interval containing }$x$,
even if it is empty and so formally not a net interval.

\begin{theorem}
\label{thm:comparable} Let $\mu $ be a self-similar measure of finite type
and let $x\in K$. Then 
\begin{equation*}
\dim _{loc}\mu (x)=\lim_{n\rightarrow \infty }\frac{\log M_{n}(x)}{n\log
\varrho },
\end{equation*}%
provided the limit exists. The lower and upper local dimensions of $\mu $ at 
$x$ can be expressed similarly in terms of $\liminf $ and $\limsup $.
\end{theorem}

\begin{proof}
Assume, first, that%
\begin{equation*}
\dim _{loc}\mu (x)=\lim_{n\rightarrow \infty }\frac{\log \mu \lbrack
x-\varrho ^{n},x+\varrho ^{n}]}{n\log \varrho }=D
\end{equation*}%
exists.

By the finite type assumption, there are constants $0<c<C$ such that $%
c\varrho ^{n}<\ell (\Delta _{n})<C\varrho ^{n}$ for all $\Delta _{n}\in 
\mathcal{F}_{n}$. Pick $j$ and $k$ such that $\varrho ^{j}<c$ and $%
2C<\varrho ^{-k}$.

If $x$ is a boundary point, then for sufficiently large $n$, $x$ is an
endpoint of $\Delta _{n}(x)$ and%
\begin{equation*}
\lbrack x-\varrho ^{n+j},x+\varrho ^{n+j}]\subseteq \Delta _{n}(x)\cup
\Delta _{n}^{\prime }(x)\subseteq \lbrack x-\varrho ^{n-k},x+\varrho ^{n-k}],
\end{equation*}%
where the notation is as in (\ref{Mnbdy}). If $x$ is not a boundary point,
then 
\begin{equation*}
\lbrack x-\varrho ^{n+j},x+\varrho ^{n+j}]\subseteq \Delta _{n}^{-}(x)\cup
\Delta _{n}(x)\cup \Delta _{n}^{+}(x)\subseteq \lbrack x-\varrho
^{n-k},x+\varrho ^{n-k}].
\end{equation*}%
In either case, 
\begin{equation*}
\mu \lbrack x-\varrho ^{n+j},x+\varrho ^{n+j}]\leq M_{n}(x)\leq \mu \lbrack
x-\varrho ^{n-k},x+\varrho ^{n-k}].
\end{equation*}

This in turn implies that 
\begin{align*}
\left( \frac{n+j}{n}\right) \left( \frac{\log \mu \lbrack x-\varrho
^{n+j},x+\varrho ^{n+j}]}{(n+j)\log \varrho }\right) & \geq \frac{\log
M_{n}(x)}{n\log \varrho } \\
& \geq \left( \frac{n-k}{n}\right) \left( \frac{\log \mu \lbrack x-\varrho
^{n-k},x+\varrho ^{n-k}]}{(n-k)\log \varrho }\right) .
\end{align*}%
The limit of the left hand side and the right hand side both go to $D$,
hence the limit of the middle expression exists and is equal to $D$.

It follows similarly that if $\lim_{n\rightarrow \infty }\log M_{n}(x)/n\log
\rho $ exists, then also 
\begin{equation*}
\dim _{loc}\mu (x)=\lim_{n\rightarrow \infty }\frac{\log \mu \lbrack
x-\varrho ^{n},x+\varrho ^{n}]}{n\log \varrho }=\lim_{n\rightarrow \infty }%
\frac{\log M_{n}(x)}{n\log \varrho }.
\end{equation*}

The arguments for the lower and upper local dimensions are similar.
\end{proof}

\subsection{Periodic points}

Recall that in \cite{HHM}, $x\in K$ is called a \textit{periodic point} if $%
x $ has symbolic representation%
\begin{equation*}
\lbrack x]=(\gamma _{0},\dots ,\gamma _{J},\theta ^{-},\theta ^{-},\dots ),
\end{equation*}%
where $\theta $ is an admissible cycle (a non-trivial path with the same
first and last letter) and $\theta ^{-}$ is the path with the last letter of 
$\theta $ deleted. We refer to $\theta $ as a \textit{period} of $x$.
Boundary points are necessarily periodic and there are only countably many
periodic points. Note that a periodic point is essential if and only if it
has a period that is a path in the essential class.

If there is a choice of $\theta $ for which $T(\theta )$ is a positive
matrix, we call $x$ a \textit{positive, periodic point. }

Of course, a period for a periodic point $x$ is not unique. For example, if $%
\theta =(\theta _{1},..,\theta _{L},\theta _{1})$ is a period, then so is $%
(\theta ^{-},\theta )$ and so is $(\theta _{2},\dots,\theta _{L},\theta
_{1},\theta _{2})$. However, these different choices for the period give the
same symbolic representation for $x$. But if $x$ is a boundary point, then $%
x $ may have two different symbolic representations, one for which $%
[x|N]=\Delta _{N}(x)$ and the other having $[x|N]=\Delta _{N}^{\prime }(x)$,
and these two representations arise from (fundamentally) different periods.

The notation $sp(T)$ means the spectral radius of the matrix $T$,%
\begin{equation*}
sp(T)=\lim_{n\rightarrow \infty }\left\Vert T^{n}\right\Vert ^{1/n}.
\end{equation*}%
We note that two periods associated with the same symbolic representation
for $x$ will have the same spectral radius. This need not be the case for
periods associated with different symbolic representations.

Here is the analogue of \cite[Proposition 4.14]{HHM} when there is no
assumption of regularity.

\begin{proposition}
\label{periodic}If $x$ is a periodic point with period $\theta $, then the
local dimension exists and is given by 
\begin{equation*}
\dim_{loc}\mu (x)=\frac{\log sp(T(\theta ))}{L(\theta ^{-})\log \varrho },
\end{equation*}%
where if $x$ is a boundary point of a net interval with two different
symbolic representations given by periods $\theta $ and $\phi $, then $%
\theta $ is chosen to satisfy $sp(T(\theta ))\geq sp(T(\phi ))$.
\end{proposition}

\begin{proof}
First, suppose $x$ is a boundary periodic point with two different symbolic
representations given by periods $\theta ,\phi $ and that $sp(T(\theta
))\geq sp(T(\phi ))$. There is no loss of generality in assuming the two
periods have the same lengths $L=L(\theta ^{-})$ and pre-period path of
length $J$. Given large $n$, let $m=[(n-J)/L]$, so $x$ has symbolic
representations 
\begin{equation*}
(\gamma _{0},\gamma _{1},\dots ,\gamma _{J-1},\underbrace{\theta ^{-},\dots
,\theta ^{-}}_{m},\theta _{1},\dots ,\theta _{t})
\end{equation*}%
and 
\begin{equation*}
(\gamma _{0},\gamma _{1}^{\prime },\dots ,\gamma _{J-1}^{\prime },%
\underbrace{\phi ^{-},\dots ,\phi ^{-}}_{m},\phi _{1},\dots ,\phi _{t})
\end{equation*}%
for suitable $t\leq L$. From Proposition \ref{HHM}, 
\begin{equation*}
\mu (\Delta _{n}(x))\sim \left\Vert T(\gamma _{0},\dots ,\gamma _{J-1},%
\underbrace{\theta ^{-},\dots ,\theta ^{-}}_{m},\theta _{1},\dots ,\theta
_{t})\right\Vert
\end{equation*}%
and 
\begin{equation*}
\mu (\Delta _{n}^{\prime }(x))\sim \left\Vert T(\gamma _{0},\dots ,\gamma
_{J-1}^{\prime },\underbrace{\phi ^{-},\dots ,\phi ^{-}}_{m},\phi _{1},\dots
,\phi _{t})\right\Vert .
\end{equation*}%
Lemma \ref{compare} implies that there are positive constants $c_{j}$,
independent of $n$, such that 
\begin{eqnarray*}
\left\Vert (T(\theta ))^{m+1}\right\Vert &\leq &\left\Vert T(\underbrace{%
\theta ^{-},\dots ,\theta ^{-}}_{m},\theta _{1},\dots ,\theta
_{t})\right\Vert \left\Vert T(\theta _{t},\dots ,\theta _{L},\theta
_{1})\right\Vert \\
&\leq &c_{1}\left\Vert T(\gamma _{0},\dots ,\gamma _{J},\underbrace{\theta
^{-},\dots ,\theta ^{-}}_{m},\theta _{1},\dots ,\theta _{t})\right\Vert \leq
c_{2}\left\Vert (T(\theta ))^{m}\right\Vert ,
\end{eqnarray*}%
and consequently,%
\begin{equation*}
c_{3}\left\Vert (T(\theta ))^{m+1}\right\Vert \leq \mu (\Delta _{n}(x))\leq
c_{4}\left\Vert (T(\theta ))^{m}\right\Vert .
\end{equation*}%
Similarly, 
\begin{equation*}
c_{3}^{\prime }\left\Vert (T(\phi ))^{m+1}\right\Vert \leq \mu (\Delta
_{n}^{\prime }(x))\leq c_{4}^{\prime }\left\Vert (T(\phi ))^{m}\right\Vert ,
\end{equation*}

If $sp(T(\theta ))>sp(T(\phi ))$, then for large enough $m$, $\left\Vert
(T(\theta ))^{m}\right\Vert >\left\Vert (T(\phi ))^{m}\right\Vert $ and hence%
\begin{eqnarray*}
c_{3}\left\Vert (T(\theta ))^{m+1}\right\Vert &\leq &\mu (\Delta
_{n}(x))+\mu (\Delta _{n}^{\prime }(x)) \\
&=&M_{n}(x)\leq 2c_{4}\left\Vert (T(\theta ))^{m}\right\Vert
\end{eqnarray*}%
Since $\left\Vert (T(\theta ))^{m}\right\Vert ^{1/m}\rightarrow sp(T(\theta
)),$ Theorem \ref{thm:comparable} gives 
\begin{equation*}
\dim _{loc}\mu (x)=\lim_{n}\frac{\log M_{n}(x)}{n\log \rho }=\lim_{m}\frac{%
\log \left\Vert T(\theta )^{m}\right\Vert }{mL\log \rho }=\frac{\log
sp(T(\theta ))}{L\log \rho }.
\end{equation*}

If, instead, $sp(T(\theta ))=sp(T(\phi ))$, then for each $n$, 
\begin{equation*}
C_{1}\max \left( \left\Vert (T(\theta ))^{m+1}\right\Vert ,\left\Vert
(T(\phi ))^{m+1}\right\Vert \right) \leq M_{n}(x)\leq C_{2}\max \left(
\left\Vert (T(\theta ))^{m}\right\Vert ,\left\Vert (T(\phi ))^{m}\right\Vert
\right) .
\end{equation*}%
As both%
\begin{equation*}
\frac{\log \left\Vert (T(\theta ))^{m}\right\Vert }{mL\log \rho },\frac{\log
\left\Vert T(\phi ))^{m}\right\Vert }{mL\log \rho }\longrightarrow
_{m\rightarrow \infty }\frac{\log sp(T(\theta ))}{L\log \rho },
\end{equation*}%
the result again follows.

If $x$ is a boundary periodic point with only one symbolic representation,
then $\Delta _{n}^{\prime }(x)$ is empty for large $n$ and the arguments are
similar, but easier.

Now, assume $x$ is not a boundary point. Then there is no loss of generality
in assuming $[x]=(\gamma _{0},\dots ,\gamma _{J},\theta ^{-},\theta
^{-},\dots )$ where the net interval $(\gamma _{0},\dots ,\gamma _{J},\theta
^{-})$ is in the interior of the net interval $(\gamma _{0},\dots ,\gamma
_{J},\theta _{1}),$ and $\theta _{1}$ is the first letter of $\theta $. This
ensures that $(\gamma _{0},\dots ,\gamma _{J},\theta _{1})$ is a common
ancestor of the two adjacent intervals to%
\begin{equation*}
(\gamma _{0},\dots ,\gamma _{J},\theta _{1},\dots ,\theta _{L},\theta
_{1},\dots ,\theta _{t})\text{ for }0\leq t\leq L=L(\theta ^{-})\text{,}
\end{equation*}%
that is at at most $2L$ levels back. More generally, $(\gamma _{0},\dots
,\gamma _{J},\underbrace{\theta ^{-},\dots ,\theta ^{-}}_{m},\theta _{1})$
is a common ancestor of the two adjacent intervals to the net intervals 
\begin{equation*}
(\gamma _{0},\dots ,\gamma _{J},\underbrace{\theta ^{-},\dots ,\theta ^{-}}%
_{m},\theta _{1},\dots ,\theta _{L},\theta _{1},\dots ,\theta _{t})\text{
for }0\leq t\leq L\text{,}
\end{equation*}%
again at most $2L$ levels back. Thus, if 
\begin{equation*}
\lbrack x|n]=(\gamma _{0},\dots ,\gamma _{J},\underbrace{\theta ^{-},\dots
,\theta ^{-}}_{m},\theta _{1},\dots ,\theta _{L},\theta _{1},\dots ,\theta
_{t}),
\end{equation*}%
then $M_{n}(x)$ is comparable to $\mu (\Delta _{n}(x)),$ with constants of
comparability independent of $n$. As above, there are $c_{j}>0$ so that%
\begin{equation*}
c_{1}\left\Vert (T(\theta ))^{m+1}\right\Vert \leq \mu (\Delta _{n}(x))\leq
c_{2}\left\Vert (T(\theta ))^{m}\right\Vert
\end{equation*}%
and the argument is completed in a similar fashion to before.
\end{proof}

\section{\protect\bigskip Local Dimensions at Truly Essential Points}

In this section we will obtain our main theoretical results on the structure
of local dimensions, analogues of those found in Section 5 of \cite{HHM}.
Because local dimensions may depend on adjacent net intervals, $\Delta
_{n}^{-}(x)$ and $\Delta _{n}^{+}(x),$ rather than only on $\Delta _{n}(x)$,
we introduce a subset of the essential class that we call the truly
essential class. We will see that this subset has full $\mu $ and Hausdorff $%
s$-measure for $s=\dim _{H}K$. Our main results state that under a weak
technical assumption the local dimensions at periodic points are dense in
the set of (upper and lower) local dimensions at truly essential points and
that the set of local dimensions at truly essential points is a closed
interval. Furthermore, we prove that there is always a truly essential point
at which the local dimension agrees with the Hausdorff dimension of the
self-similar set and we give criteria for when the measure $\mu $ is
absolutely continuous with respect to the Hausdorff measure.

\label{sec:TEP}

\subsection{Truly essential points}

\begin{definition}
Suppose $K$ is the self-similar set associated with an IFS of finite type.

\begin{enumerate}
\item We will say that $x\in K$ is a \textbf{boundary essential point} if $x$
is a boundary point of $\Delta _{n}(x)\in \mathcal{F}_{n}$ for some $n,$ and
both $\Delta _{n}(x)$ and the other $n$'th level net interval containing $x$%
, $\Delta _{n}^{\prime }(x),$ are essential (where if $\Delta _{n}^{\prime
}(x)$ is empty we understand it to be essential).

\item We will say that $x\in K$ is an \textbf{interior essential point} if $%
x $ is not a boundary point and there exists an essential net interval with $%
x$ in its interior.

\item We call $x$ a \textbf{truly essential point} if it is either an
interior essential point or a boundary essential point.
\end{enumerate}
\end{definition}

Obviously, truly essential points are essential and if $x$ is in the
interior of some essential net interval, then it is truly essential. In
particular, any essential point that is not truly essential must be a
boundary point. Hence there can be only countably many of these and they are
periodic.

Any point in the relative interior of the essential class (with respect to
the space $K$), is either contained in the interior of some essential
interval, or is a boundary essential point. Hence the relative interior of
the essential class is equal to the set of truly essential points. If the
essential class is a (relatively) open set, then the essential class
coincides with the truly essential class. This is the situation, for
example, with the Bernoulli convolutions and Cantor-like measures discussed
in Sections \ref{sec:Pisot} and \ref{sec:Cantor}. Another IFS where the set
of essential points is equal to the set of truly essential points is given
in Example \ref{Ex:notpos}.

However, as the example below demonstrates, these two sets need not be equal.

\begin{example}
\label{ex:new} Consider the maps $S_{i}(x)=x/4+d_{i}/8$ with $d_{i}=i$ for $%
i=0,\dots ,3$, $d_{4}=5$, and $d_{5}=6$. The reduced transition diagram has
4 reduced characteristic vectors. The reduced characteristic vectors are:

\begin{itemize}
\item Reduced characteristic vector 1: $(1,(0))$

\item Reduced characteristic vector 2: $(1/2,(0))$

\item Reduced characteristic vector 3: $(1/2,(0,1/2))$

\item Reduced characteristic vector 4: $(1/2, (1/2))$
\end{itemize}

The transition maps are:

\begin{itemize}
\item RCV $1 \to [2a, 3a, 3b, 3c, 4a, 2rb, 3d, 4b]$

\item RCV $2 \to [2, 3a, 3b, 3c]$

\item RCV $3 \to [3a, 3b, 3c, 3d]$

\item RCV $4 \to [4a, 2, 3, 4b]$
\end{itemize}

By this we mean, for example, that the reduced characteristic vector 1 has 8
children. Listed in order from left to right, they are the reduced
characteristic vectors $1,2,3,3,3,4,2,3,4$ etc. By $3a$ we mean the first
occurance of the child of type $3$, $3b$ the second, etc. If there is only
one child of that type we do not need to distinquish them. It is possible
that the transitions matrices are different for different children of the
same type. These also help to distinquish paths unamibiguously. See Example %
\ref{Ex:notpos}. It is easy to see from the
transition maps that the essential class is $\{3a,3b,3c,3d\}$. See Figure %
\ref{fig:Pic7} for the transition diagram. 
\begin{figure}[tbp]
\includegraphics[scale=0.5]{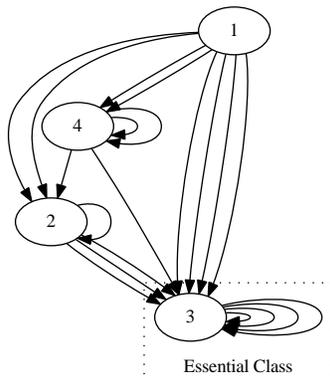}
\caption{Transition diagram for Example \protect\ref{ex:new}}
\label{fig:Pic7}
\end{figure}

Consider the boundary periodic point $x$ with symbolic representation $%
(1,4a,4a,4a,\dots )$, $4a$ being the left most child of $4$. This also has
symbolic representation $(1,3c,3d,3d,\dots ),$ $3d$ being the right-most
child of $3$. One of these symbolic representations is in the essential
class, whereas the other is not. As such, this point is an essential point,
but it is not a truly essential point.
\end{example}

The significance of an interior essential point $x$ is that $\Delta
_{n}^{-}(x),\Delta _{n}(x)$ and $\Delta _{n}^{+}(x)$ have a common essential
ancestor for some $n$. Conversely, if $\Delta _{n}^{-}(x),\Delta _{n}(x)$
and $\Delta _{n}^{+}(x)$ have a common essential ancestor for some $n$, then 
$x$ belongs to the relative interior of the essential class and thus is
truly essential.

A periodic point $x$ that is an interior essential point admits a period $%
\theta $ with the property that if $[x]=(\gamma ,\theta ^{-},\theta
^{-},\dots )$ and $\theta _{1}$ is the first letter of $\theta ,$ then the
net interval (with symbolic representation) $(\gamma ,\theta )$ is in the
interior of the net interval $(\gamma ,\theta _{1})$. We will call such a
period $\theta $ \textit{truly essential}. Equivalently, $\theta $ is truly
essential if and only if $\theta $ is a path that does not consist solely of
right-most descendents or solely of left-most descendents.

It was shown \cite[Proposition 4.5]{HHM} that under the assumption that the
self-similar set was an interval, the essential class had full Lebesgue
measure. In fact, this is true for the truly essential class, with Lebesgue
measure replaced by either the self-similar measure $\mu $ or the
(normalized) Hausdorff $s$-measure, where $s=\dim _{H}K$, as the next
Proposition shows. To prove this, we first need some preliminary lemmas.

\begin{lemma}
There exists an integer $J$ such that for each net interval $\Delta \in 
\mathcal{F}_n$ there exists a $\sigma \in \mathcal{A}^{J+n}$ with $%
\sigma([0,1]) \subseteq \Delta$.
\end{lemma}

\begin{proof}
Consider a net interval $\Delta \in \mathcal{F}_{n}$. As there is some $x\in
K$ in the interior of $\Delta $, there is an index $t$ and $\sigma \in 
\mathcal{A}^{n+t}$ such that $\sigma ([0,1])\subseteq \Delta $. Since $%
\sigma (0)$ is not isolated in $K$, there must be a level $n+T$ net interval 
$\Delta _{0}\subseteq \Delta ,$ with left end $S_{\sigma }(0)$. Choose the
index $T$ minimal with this property.

Assume $\Delta $ has symbolic representation $(\gamma _{0},\gamma
_{1},\dots,\gamma _{n})$ and $\Delta _{0}$ has representation $(\gamma
_{0},\gamma _{1},\dots,\gamma _{n},\dots,\gamma _{n+T})$. Let $\Delta
^{\prime }$ $=$ $(\gamma _{0},\chi _{1},\dots,\chi _{m-1},\gamma _{n})\in 
\mathcal{F}_{m}$ be any other net interval with symbolic representation
ending with the same characteristic vector $\gamma _{n}$. It will also have
a descendent, $\Delta _{0}^{\prime }$, with representation ending with the
path $(\gamma _{n},\dots,\gamma _{n+T})$. Since $0$ is in the neighbour set
of $\Delta _{0}$, the same is true for $\Delta _{0}^{\prime }$ and thus its
left endpoint is an image of $0$ under $S_{\tau }$ for some $\tau \in 
\mathcal{A}^{m+T}$. As the pairs $(\Delta ,$ $\Delta ^{\prime })$ and $%
(\Delta _{0},\Delta _{0}^{\prime })$ have the same finite type structure (up
to normalization), it follows that $S_{\tau }([0,1])\subseteq \Delta
^{\prime }$. Hence $\Delta ^{\prime }$ has the same minimal index $T$, in
other words, $T$ depends only upon the final characteristic vector
associated with $\Delta $.

As there are only finitely many characteristic vectors, we can take $J$ to
be the maximum of these indices $T$ taken over all the characteristic
vectors.
\end{proof}

\begin{lemma}
\label{lem:H(Delta)} There exists a positive constant $c$ such that for all $%
\Delta \in \mathcal{F}_{n}$ and all $n$ we have $c\varrho ^{sn}\leq
H^{s}(\Delta \cap K)\leq \varrho ^{sn}.$
\end{lemma}

\begin{proof}
Fix an $n$'th level net interval $\Delta \in \mathcal{F}_{n}$. By
construction, there exists some $\sigma \in \mathcal{A}^{n}$ such that $%
\Delta \subseteq \sigma ([0,1])$. Then $\Delta \cap K\subseteq \sigma (K)$
and hence $H^{s}(\Delta \cap K)\leq \varrho ^{ns}H^{s}(K)=\varrho ^{ns}$.

Choose $J$ as in the previous lemma. Then there exists some $\tau \in 
\mathcal{A}^{n+J}$ such that $\tau ([0,1])\subseteq \Delta $. Hence $\tau
(K)\subseteq \Delta \cap K$ and therefore $\varrho ^{(J+n)s}H^{s}(K)\leq
H^{s}(\Delta \cap K)$. Taking $c=\varrho ^{Js}>0,$ we are done.
\end{proof}

\begin{corollary}
We have $\dim _{H}(\Delta \cap K)=s$ for all net intervals $\Delta $.
\end{corollary}

\begin{proof}
This is immediate since $0<H^{s}(\Delta \cap K)<\infty $.
\end{proof}

\begin{proposition}
\label{SizeofTE}Suppose $\mu $ is a self-similar measure of finite type,
with support $K$ of Hausdorff dimension $s$. The set of points in $K$ that
are not truly essential is a subset of a closed set having zero $\mu $ and $%
H^{s}$-measure.
\end{proposition}

\begin{proof}
As we already observed, every net interval contains a descendent net
subinterval that is essential. This essential net interval contains some $%
x\in K$ in its interior and hence contains a further subinterval that is in
its interior. For the purposes of this proof, we will call this an interior
essential net interval. The finite type property ensures we can always find
an interior essential net subinterval within a bounded number of
generations, say at most $J$.

We claim that there exists some $\lambda >0$ such that the proportion of the
measure of this net subinterval to the measure of the original interval is $%
\geq \lambda $. This is because all $J$'th level descendent net subintervals
have comparable measure to the original net interval. For measure $H^{s}$,
this property is shown in Lemma \ref{lem:H(Delta)} and for the measure $\mu $
it follows from the definition.

We now exhibit a Cantor-like construction. We begin with $[0,1]$. Consider
the first level at which there is a net interval that is interior essential.
Remove the interiors of all the net intervals of this level that are
interior essential. The resulting closed subset of $[0,1]$ is a finite union
of closed intervals, say $C_{1}$, whose measures, either $H^{s}$ or $\mu $,
total at most $1-\lambda $. We repeat the process of removing the interiors
of the interior essential net intervals at the next level at which there are
interior essential, net intervals in each of the intervals of $C_{1}$. The
resulting closed subset now has measure at most $(1-\lambda )^{2}$.

After repeating this procedure $k$ times one can see that the non-interior
essential points are contained in a finite union of closed intervals,
denoted $C_{k},$ whose total measure is at most $(1-\lambda )^{k}$. It
follows that the non-interior essential points are contained in the closed
set $\bigcap\limits_{k=1}^{\infty }C_{k},$ and this set has both $\mu $ and $%
H^{s}$-measure $0$.
\end{proof}

\begin{remark}
Observe that we have actually proven that the complement of the interior of
the essential class (in $K$) has $\mu $ and $H^{s}$-measure zero.
\end{remark}

Another consequence of Lemma \ref{lem:H(Delta)} is to obtain a new formula
for the Hausdorff dimension of self-similar set of finite type. In \cite{NW}
a formula was given which required knowing the complete transition graph. In
fact, it suffices to know the transition graph of the essential
characteristic vectors. For the purpose of this proof we introduce the
following notation: Let $\gamma _{1},\dots ,\gamma _{r}$ be a complete list
of the reduced characteristic vectors. Define a $r\times r$ matrix $I$ by $%
(I)_{jk}=$ the number of children of $\gamma _{j}$ that are of type $\gamma
_{k}$. We call $I$ the \textit{incidence matrix of the essential class}.

\begin{proposition}
\label{dimH}Let $K$ be a self-similar set of finite type and let $I$ be the
incidence matrix of the essential class. Then 
\begin{equation*}
s=\dim _{H}K=\frac{\log (sp(I))}{\left\vert \log \varrho \right\vert }.
\end{equation*}
\end{proposition}

\begin{example}
Consider the example $S_{j}(x)=x/3+b_{j}$ with $b_{j}\in \{0,2/87,2/3\}$.
This IFS has 2280 reduced characteristic vectors, hence to compute the
dimension using the full set of reduced characteristic vectors would require
finding the eigenvalues of a $2280\times 2280$ matrix. But there are only 2
essential vectors and the incidence matrix of the essential class is equal
to 
\begin{equation*}
\left[ 
\begin{matrix}
2 & 1 \\ 
1 & 2%
\end{matrix}%
\right] .
\end{equation*}%
Using the proposition above one can easily deduce that the dimension of the
self-similar set is $1$, although the set is not the full interval $[0,1]$.
\end{example}

\begin{proof}[Proof of Proposition \protect\ref{dimH}.]
Choose $\Delta _{0}$ an essential net interval of level $n_{0}$ with the
property that all essential characteristic vectors are descendents of $%
\Delta _{0}$ at level $N+n_{0}$. It can be seen from the proof of Lemma 6.4
of \cite{F2} that such a net interval exists. For $n\geq N,$ let 
\begin{equation*}
E_{n}=\{\Delta \in \mathcal{F}_{n+n_{0}}:\Delta \subseteq \Delta _{0}\}.
\end{equation*}

From Lemma \ref{lem:H(Delta)} we have (for $\left\vert E_{n}\right\vert $
denoting the cardinality of $E_{n})$,%
\begin{equation*}
c\varrho ^{sn}\left\vert E_{n}\right\vert \leq \sum_{\Delta \in
E_{n}}H^{s}(\Delta \cap K)\leq C\varrho ^{sn}\left\vert E_{n}\right\vert
\end{equation*}%
for positive constants $c,C$. Since the sets $\Delta \in E_{n}$ have
disjoint interiors, 
\begin{equation*}
0<H^{s}(\Delta _{0}\cap K)=\sum_{\Delta \in E_{n}}H^{s}(\Delta \cap
K)<\infty ,
\end{equation*}%
thus there are positive constants $A,B$ such that%
\begin{equation*}
A\leq \varrho ^{sn}\left\vert E_{n}\right\vert \leq B\text{ for all }n.
\end{equation*}%
Consequently,%
\begin{equation*}
\frac{\frac{1}{n}\log \left\vert E_{n}\right\vert }{\left\vert \log \varrho
\right\vert }+\frac{\log A}{n\log \varrho }\geq s\geq \frac{\frac{1}{n}\log
\left\vert E_{n}\right\vert }{\left\vert \log \varrho \right\vert }+\frac{%
\log B}{n\log \varrho }.
\end{equation*}

Without loss of generality we can assume $\Delta _{0}$ has symbolic
representation with last letter $\gamma _{1}$. Then $\left\vert
E_{n}\right\vert $ is the sum of the entries of row $1$ of $I^{n}$, so 
\begin{equation*}
\left\vert E_{n}\right\vert =\left\Vert [1,0,\dots ,0]I^{n}\right\Vert
=\left\Vert [1,0,\dots ,0]I^{N}I^{n-N}\right\Vert .
\end{equation*}%
But $[1,0,\dots ,0]I^{N}$ is a vector with all non-zero entries since $%
\Delta _{0}$ has all the essential characteristic vectors as descendents at
level $N+n_{0}$. Hence $\left\vert E_{n}\right\vert \sim \left\Vert
I^{n-N}\right\Vert $ and since 
\begin{equation*}
\frac{\frac{1}{n}\log \left\Vert I^{n-N}\right\Vert }{\left\vert \log
\varrho \right\vert }\rightarrow \frac{\log (sp(I))}{\left\vert \log \varrho
\right\vert },
\end{equation*}%
we deduce that this is the value of $s$.
\end{proof}

\subsection{Positive row property}

Throughout the remainder of this section, we will assume, without loss of
generality, that $\Delta \subsetneqq \widehat{\Delta }$ whenever the net
interval $\Delta $ is a child of $\widehat{\Delta }$. To see that this
assumption is without loss of generality, we note that as $\mu $ is of
finite type, there will be an integer $N$ such that all net intervals will
have at least two descendents $N$ levels deeper. Consider the new IFS with
contractions $S_{i_{1}}\circ \dots \circ S_{i_{N}}$ and probabilities $%
p_{i_{1}}\dots p_{i_{N}}$. This IFS gives rise to the same self-similar
measure $\mu $. Moreover, the set of net intervals of level $kN$ of the
original construction are precisely the level $k$ net intervals in the new
construction. This new construction has the desired property.

In the theorems of this subsection we will also assume that the self-similar
measure of finite type has the property that each essential primitive
transition matrix has a non-zero entry in each row. This is the weak
technical condition referred to in the introduction and we call it the 
\textit{positive row property}. The property holds automatically when the
self-similar set $K=[0,1]$ (see \cite[Sec. 3.2]{HHM}), such as for (even
non-regular) Bernoulli convolutions and Cantor-like measures. This stronger
assumption is not necessary, though, as we see in Example \ref{Ex:pos}.

The positive row property can fail to hold when $K\neq \lbrack 0,1]$ and can
even fail when there is a positive essential transition matrix, as the
example below demonstrates.

\begin{example}
\label{Ex:notpos} Consider the self-similar measure associated with the IFS 
\begin{equation*}
\{S_{j}(x)=x/3+d_{j}:d_{j}=0,4/9,5/9,2/3\}
\end{equation*}%
and uniform probabilities. This measure is of finite type. Its support is a
proper subset of $[0,1]$ since $(1/3,4/9)\cap K$ is empty. However, if $I=$ $%
[2/3,1],$ then $I\subseteq \bigcup_{j=0}^{3}S_{j}(I)$ and this implies that $%
[2/3,1]\subseteq K$. Thus $K$ has positive Lebesgue measure. The reduced
characteristic vectors are:

\begin{itemize}
\item Reduced characteristic vector 1: $(1, (0))$

\item Reduced characteristic vector 2: $(1/3, (0))$

\item Reduced characteristic vector 3: $(1/3, (0, 1/3))$

\item Reduced characteristic vector 4: $(1/3, (0, 1/3, 2/3))$

\item Reduced characteristic vector 5: $(1/3, (1/3, 2/3))$

\item Reduced characteristic vector 6: $(1/3, (2/3))$
\end{itemize}

The transition maps are:

\begin{itemize}
\item RCV $1 \to [1, X, 2, 3, 4, 5, 6]$

\item RCV $2 \to [1]$

\item RCV $3 \to [2, 3, 4]$

\item RCV $4 \to [4, 4, 4]$

\item RCV $5 \to [4, 4, 4]$

\item RCV $6 \to [4, 5, 6]$
\end{itemize}

The `$X$' denotes that between the child of type $1$ and the child of type $%
2,$ in the parent $1$, there is an interval $[h_{j},h_{j+1}]$ that is not a
net interval, as $(h_{j},h_{j+1})\cap K=\emptyset $.

The transition diagram is shown in Figure \ref{fig:zerorow}. There are three
(non-reduced) essential characteristic vectors denoted $4a,4b,4c$ and one
reduced characteristic vector $4$. The primitive transition matrices for the
essential class are given below. For $x=a,b,c$, the matrix $T(4,4x)$ is any
of $T(4a,4x),T(4b,4x)$ or $T(4c,4x)$ as these three matrices coincide. Note
that $T(4,4a)$ has a row of zeroes, while the essential transition matrix $%
(T(4c,4b)T(4b,4b)T(4b,4c))^{2}$ is positive. Hence this example does not
satisfy the positive row property, although there is a positive essential
transition matrix.

\begin{equation*}
T(4,4a)=\left[ 
\begin{array}{ccc}
\frac{1}{4} & \frac{1}{4} & \frac{1}{4} \\ 
0 & 0 & 0 \\ 
0 & 0 & \frac{1}{4}%
\end{array}%
\right] ,T(4,4b)=\left[ 
\begin{array}{ccc}
\frac{1}{4} & \frac{1}{4} & 0 \\ 
0 & 0 & \frac{1}{4} \\ 
0 & \frac{1}{4} & 0%
\end{array}%
\right] ,T(4,4c)=\left[ 
\begin{array}{ccc}
\frac{1}{4} & 0 & 0 \\ 
0 & \frac{1}{4} & \frac{1}{4} \\ 
\frac{1}{4} & 0 & 0%
\end{array}%
\right]
\end{equation*}

\begin{figure}[tbp]
\includegraphics[scale=0.7]{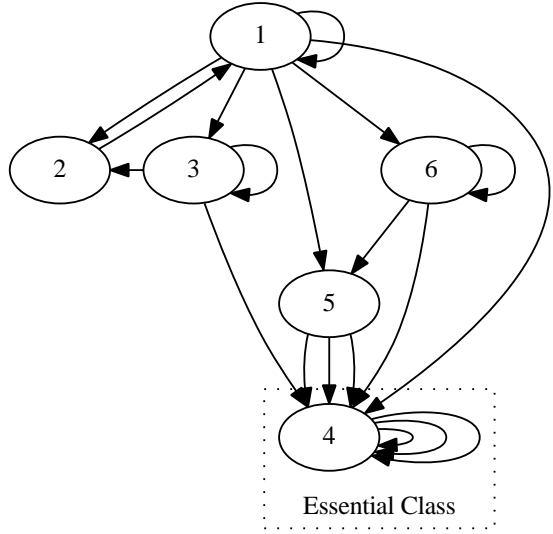}
\caption{Transition diagram for Example \protect\ref{Ex:notpos}}
\label{fig:zerorow}
\end{figure}

If $x$ is an essential point, but not a truly essential point, then $x$
cannot be in the interior of an essential interval. This means $x$ will be
on the boundary of both an essential net interval, $\Delta _{n}(x)$, and a
non-essential net interval, $\Delta _{n}^{\prime }(x),$ for all $n$
sufficiently large. It is easy to see from the transition maps that there
are no such points. Hence the truly essential set coincides with the
essential set.
\end{example}

\subsection{Main Results}

We begin by establishing the existence of special paths that we call truly
essential.

\begin{lemma}
Suppose $\mu $ is a self-similar measure of finite type satisfying the
positive row property. Given any two essential characteristic vectors, $%
\gamma _{1},\gamma _{2}$, there is a positive, essential path $\eta =(\eta
_{1},\dots ,\eta _{k}),$ that does not consist of solely left-most
descendents, or solely right-most descendents, and having $\eta _{1}=$ $%
\gamma _{1}$ and $\eta _{k}=\gamma _{2}$.
\end{lemma}

\begin{proof}
In Proposition 4.12 of \cite{HHM} it is shown that there is an admissible
essential path $\eta _{0}$ that begins and ends at $\gamma _{1}$ and is
positive. Since any net interval contains an element of $K$ in its interior,
there must be an essential path, $\eta ^{\prime },$ beginning with $\gamma
_{1}$ and ending at, say $\chi $, which does not consist of solely left-most
or solely right-most descendents. Now take any essential path $\eta ^{\prime
\prime }$ from $\chi $ to $\gamma _{2}$. Put $\eta =\eta _{0}\eta ^{\prime
}\eta ^{\prime \prime }$. This is a positive path since the product (in
either order) of any positive matrix by a matrix with a non-zero entry in
each row and column is again positive.
\end{proof}

We will call a path $\eta =(\eta _{1},\dots ,,\eta _{s}),$ as described in
the lemma above, a \textit{truly essential, positive path}. If $\Delta
=(\gamma _{0},\dots ,\gamma _{N},\eta _{1},\dots ,\eta _{s})$, then the two
adjacent intervals of $\Delta $ are both descendents of the essential
interval $(\gamma _{0},\dots ,\gamma _{N},\eta _{1})$. Consequently any $x$
whose symbolic representation begins $(\gamma _{0},\dots ,\gamma _{N},\eta )$
is truly essential. Further, if $x$ is a periodic point with period $\theta
\eta \phi $ for some $\theta $ and $\phi ,$ then $x$ is an interior
essential point.

\textbf{Notation:} For the remainder of this section, $F$ will denote a
fixed, finite set of truly essential, positive paths with the property that
given any two essential characteristic vectors, there is a path in $F$
joining them in either order.

\begin{theorem}
Suppose $\mu $ is a self-similar measure of finite type satisfying the
positive row property. Then the set of lower local dimensions of $\mu $ at
interior essential, positive, periodic points is dense in the set of all
local dimensions of $\mu $ at truly essential points. A similar statement
holds for the (upper) local dimensions.
\end{theorem}

\begin{proof}
We will first assume $\dim_{loc}\mu (x)$ exists. The arguments for upper and
lower local dimensions are similar.

\textbf{Step 1:} To begin, we will show that if $x$ is a boundary essential
(necessarily periodic) point, then its local dimension can be approximated
by that of a truly essential, positive, periodic point.

First, suppose $x$ has two different symbolic representations, say 
\begin{equation*}
(\gamma _{0},\gamma _{1},\dots,\gamma _{J},\theta ^{-},\theta ^{-},\dots )%
\text{ and }(\gamma _{0},\gamma _{1}^{\prime },\dots,\gamma _{J^{\prime
}}^{\prime },\phi ^{-},\phi ^{-},\dots ).
\end{equation*}%
%
%
There is no loss of generality in assuming the two periods have the same
lengths $L=L(\theta ^{-})$ and pre-period path of length $J$. 
Without loss of generality assume that $sp(T(\theta ))\geq
sp(T(\phi )),$ so Proposition \ref{periodic} gives 
\begin{equation*}
\dim_{loc}\mu (x)=\frac{\log sp(T(\theta ))}{L(\theta ^{-})\log \rho }.
\end{equation*}%
Let $\eta \in F$ be a truly essential, positive path chosen so that $\theta
^{-}\eta \theta $ is well defined. Consider the periodic point 
\begin{equation*}
\lbrack y_{n}]=(\gamma _{0},\dots ,\gamma _{J},\psi _{n}^{-},\psi
_{n}^{-},\dots ),
\end{equation*}%
where $\psi _{n}=\underbrace{({\theta ^{-},\dots, \theta ^{-}}}_{n}\eta
,\theta _{1}),$ $\theta _{1}$ being the first letter of $\theta $. As noted
above, this construction produces an interior essential point and therefore 
\begin{equation*}
\dim_{loc}\mu (y_{n})=\frac{\log sp(T(\psi _{n}))}{L(\psi _{n}^{-})\log \rho 
}.
\end{equation*}

As $\eta $ is a positive path, Lemma \ref{compare} implies%
\begin{equation*}
sp(T(\psi _{n}))\leq \left\Vert T(\psi _{n})\right\Vert \leq c_{1}\left\Vert
(T(\theta ))^{n}\right\Vert \leq c_{2}sp(T(\theta ^{n}))
\end{equation*}%
and%
\begin{equation*}
sp(T(\psi _{n}))\geq c_{3}\left\Vert T(\psi _{n})\right\Vert \geq
c_{4}sp(T(\theta ))^{n},
\end{equation*}%
where the constants are positive and independent of $n$. It follows that 
\begin{equation*}
\dim_{loc} \mu (y_{n})=\frac{\log C_{n}sp(T(\theta ))^{n}}{(nL(\theta
^{-})+L(\eta ))\log \rho }
\end{equation*}%
where the constants, $C_{n},$ are bounded above and bounded below from zero. 
Hence

\begin{eqnarray*}
\dim_{loc} \mu (y_{n}) &=&\frac{\log C_{n}}{(nL(\theta ^{-})+L(\eta ))\log
\rho }+\frac{\log sp(T(\theta ))}{(L(\theta ^{-})+\frac{1}{n}L(\eta ))\log
\rho } \\
&=&\frac{\log C_{n}}{(nL(\theta ^{-})+L(\eta ))\log \rho }+\dim_{loc} \mu (x)%
\frac{L(\theta ^{-})}{L(\theta ^{-})+\frac{1}{n}L(\eta )} \\
&\rightarrow &\dim_{loc} \mu (x)\text{ as }n\rightarrow \infty .
\end{eqnarray*}

The case where $x$ has a unique representation is similar. This completes
step 1.

\textbf{Step 2:} Now, suppose $x$ is an interior essential point with $%
\Delta _{N}(x)$ and its two adjacent $N$'th level intervals having common
essential ancestor at level $J$. If the symbolic representation for $x$
begins with the path $(\gamma _{0},\dots,\gamma _{J})$, then for any $n>N$,
all three of $\Delta _{n}(x)$, $\Delta _{n}^{+}(x),\Delta _{n}^{-}(x)$ have
symbolic representation also beginning with $(\gamma _{0},\dots,\gamma _{J})$%
.

Without loss of generality assume 
\begin{equation*}
\max \left\{ \mu (\Delta _{n}(x)),\mu (\Delta _{n}^{+}(x)),\mu (\Delta
_{n}^{-}(x))\right\} =\mu (\Delta _{n}^{+}(x))
\end{equation*}%
along a subsequence not renamed. (The other cases are similar.) Of course,
then we have

\begin{equation*}
\mu (\Delta _{n}^{+}(x))\leq M_{n}(x)\leq 3\mu (\Delta _{n}^{+}(x)).
\end{equation*}

Suppose 
\begin{equation*}
\Delta _{n}^{+}(x)=(\gamma _{0},\dots,\gamma _{J},\chi^{(n)}
_{J+1},\dots,\chi^{(n)} _{n})
\end{equation*}%
and let $\eta $ be a path in $F$ joining $\chi^{(n)} _{n}$ to $\gamma _{J}$.
We remark here that the $\chi_{J+1}^{(n)}, \dots, \chi_n^{(n)}$ will depend
on $n$ as $\Delta_n^{+}$ may not be a descendent of $\Delta_{n-1}^+$. Put 
\begin{equation*}
\theta _{n}=(\gamma _{J},\chi^{(n)} _{J+1},\dots,\chi^{(n)} _{n-1},\eta )
\end{equation*}%
and denote by $y_{n}$ the interior essential, positive, periodic point with
symbolic representation 
\begin{equation*}
\lbrack y_{n}]=(\gamma _{0},\dots,\gamma _{J-1},\theta _{n}^{-},\theta
_{n}^{-},\dots)\text{.}
\end{equation*}%
Of course,%
\begin{equation*}
\dim_{loc}\mu (y_{n})=\frac{\log sp(T(\theta _{n}))}{L(\theta _{n})\log \rho 
}.
\end{equation*}

Lemma \ref{compare} implies there is a constant $c_{1}>0,$ independent of $n$%
, such that%
\begin{equation*}
c_{1}\left\Vert T(\theta _{n})\right\Vert \leq sp(T(\theta _{n}))\leq
\left\Vert T(\theta _{n})\right\Vert .
\end{equation*}%
As $L(\theta _{n})=n-J+L(\eta )$, it follows by similar reasoning to the
above that%
\begin{equation}
\dim_{loc}\mu (y_{n})=\frac{C_{n}}{n\log \rho }+\frac{\log \left\Vert
T(\theta _{n})\right\Vert }{n\log \rho }\rightarrow \lim_{n}\frac{\log
\left\Vert T(\theta _{n})\right\Vert }{n\log \rho }.  \label{dimyn}
\end{equation}

Yet another application of Lemma \ref{compare} shows that 
\begin{align*}
M_{n}(x)& \leq 3\mu (\Delta _{n}^{+}(x))\leq c_{1}\left\Vert T(\gamma
_{0},\dots,\gamma _{J},\chi^{(n)} _{J+1},\dots,\chi^{(n)} _{n})\right\Vert \\
& \leq c_{2}\left\Vert T(\gamma _{J},\chi^{(n)} _{J+1},\dots,\chi^{(n)}
_{n},\eta )\right\Vert \leq c_{3}\left\Vert T(\theta _{n})\right\Vert
\end{align*}%
where the constants are independent of $n$, and similarly 
\begin{equation*}
M_{n}(x)\geq \mu (\Delta _{n}^{+}(x))\geq c\left\Vert T(\theta
_{n})\right\Vert .
\end{equation*}%
Thus $M_{n}(x)\sim \left\Vert T(\theta _{n})\right\Vert $ and hence it
follows from (\ref{dimyn}) that

\begin{equation*}
\dim_{loc}\mu (x)=\lim_{n}\frac{\log \left\Vert M_{n}(x)\right\Vert }{n\log
\rho }=\lim_{n}\frac{\log \left\Vert T(\theta _{n})\right\Vert }{n\log \rho }%
=\lim_{n}\dim_{loc}\mu (y_{n})
\end{equation*}
\end{proof}

\begin{theorem}
\label{limit}Suppose $\mu $ is a self-similar measure of finite type
satisfying the positive row property. Assume $(x_{n})$ are interior
essential, positive, periodic points. There there is an interior essential
point $x$ such that 
\begin{eqnarray*}
\overline{\dim }_{loc}\mu (x) &=&\limsup_{n}\dim _{loc}\mu (x_{n}) \\
\underline{\dim }_{loc}\mu (x) &=&\liminf_{n}\dim _{loc}\mu (x_{n})\text{.}
\end{eqnarray*}
\end{theorem}

\begin{proof}
This is similar to the proof of Theorem 5.5 of \cite{HHM} with some
technical complications that we highlight here. To begin, suppose $x_{n}$
has truly essential, positive period $\theta _{n}$ where, without loss of
generality, 
\begin{equation*}
\left\vert \frac{\log sp(T(\theta _{2n}))}{L(\theta _{2n}^{-})}-\limsup_{k}%
\frac{\log sp(T(\theta _{k}))}{L(\theta _{k}^{-})}\right\vert <\frac{1}{n},
\end{equation*}

\begin{equation*}
\left\vert \frac{\log sp(T(\theta _{2n+1}))}{L(\theta _{2n+1}^{-})}%
-\liminf_{k}\frac{\log sp(T(\theta _{k}))}{L(\theta _{k}^{-})}\right\vert <%
\frac{1}{n},
\end{equation*}%
all even labelled paths $\theta _{2n}^{-}$ have the same first letter and
the same last letter, and similarly for the odd labelled paths.

Choose truly essential, positive paths $\eta ^{e}$ and $\eta ^{o}$, from the
finite set $F$ so that $\eta ^{o}$ joins the last letter of an odd path to
the first letter of even path and $\eta ^{e}$ does the opposite.

Let $L_{n}=2L(\theta _{n+1}^{-})+L(\eta ^{o})+L(\eta ^{e})$. Choose $C_{n}$
such that for all $\ell \leq L_{n}$ and $j>\ell $ we have $\mu (\Delta
_{j})\geq C_{n}\mu (\Delta _{j-\ell })$ when $\Delta _{j-\ell }\in \mathcal{F%
}_{j-\ell }$ is the ancestor of $\Delta _{j}\in \mathcal{F}_{j}$. Now choose 
$k_{n}\geq 2$ sufficiently large so that in addition to the requirements of $%
k_{n}$ in the proof of \cite[Thm. 5.5]{HHM}, we also have%
\begin{equation*}
\text{ }\frac{3+\log C_{n}}{k_{n}}\rightarrow 0\text{.}
\end{equation*}%
Suppose $x\in K$ has symbolic representation 
\begin{equation*}
\lbrack x]=(\gamma _{0},\dots,\gamma _{J},\underbrace{\theta
_{1}^{-},\dots,\theta _{1}^{-}}_{k_{1}},\eta ^{o},\underbrace{\theta
_{2}^{-},\dots,\theta _{2}^{-}}_{k_{2}},\eta ^{e},\dots).
\end{equation*}%
We remark that as $\eta ^{e}$ and $\eta ^{o}$ are truly essential paths, the
point $x$ is interior essential.

Suppose 
\begin{equation*}
\lbrack x|j]=(\gamma _{0},\dots,\gamma _{J},\underbrace{\theta
_{1}^{-},\dots,\theta _{1}^{-}}_{k_{1}},\eta ^{o},\dots,\underbrace{\theta
_{n}^{-},\dots,\theta _{n}^{-}}_{k_{n}},\eta ^{\prime },(\theta
_{n+1}^{-},\theta _{n+1}^{-}))
\end{equation*}%
where $\eta ^{\prime }$ is either $\eta ^{o}$or $\eta ^{e},$ as appropriate,
and the notation $(\theta _{n+1}^{-},\theta _{n+1}^{-})$ means any subpath
of the path $\theta _{n+1}^{-},\theta _{n+1}^{-}$. As $\eta ^{\prime }$ is
truly essential, $\Delta _{j}(x)$ and its two adjacent intervals have common
ancestor 
\begin{equation*}
\lbrack x|j-\ell ]=(\gamma _{0},\dots,\gamma _{J},\underbrace{\theta
_{1}^{-},\dots,\theta _{1}^{-}}_{k_{1}},\eta ^{o},\dots,\underbrace{\theta
_{n}^{-},\dots,\theta _{n}^{-}}_{k_{n}})
\end{equation*}%
for some $\ell \leq 2L(\theta _{n+1}^{-})+L(\eta ^{\prime })$. If, instead, 
\begin{equation*}
\lbrack x|j]=(\gamma _{0},\dots,\gamma _{J},\underbrace{\theta
_{1}^{-},\dots,\theta _{1}^{-}}_{k_{1}},\eta ^{o},\dots,\underbrace{\theta
_{n}^{-},\dots,\theta _{n}^{-}}_{k_{n}},\eta ^{\prime },\underbrace{\theta
_{n+1}^{-},\dots,\theta _{n+1}^{-}}_{p_{n}},(\theta _{n+1}^{-},\eta ^{\prime
\prime }))
\end{equation*}%
with $2\leq p_{n}\leq k_{n+1}$ (where we may include a subset of $\eta
^{\prime \prime }=\eta ^{e}$ or $\eta ^{o}$ if $p_{n}=k_{n+1}-1$), then $%
\Delta _{j}(x)$ and its two adjacent intervals have common ancestor 
\begin{equation*}
\lbrack x|j-\ell ]=(\gamma _{0},\dots,\gamma _{J},\underbrace{\theta
_{1}^{-},\dots,\theta _{1}^{-}}_{k_{1}},\eta ^{o},\dots,\eta ^{\prime },%
\underbrace{\theta _{n+1}^{-},\dots,\theta _{n+1}^{-}}_{p_{n}-1})
\end{equation*}%
where $\ell \leq 2L(\theta _{n+1}^{-})+L(\eta ^{\prime \prime })$.

In either case, for all such $j$, there is some $\ell \leq L_{n}$ such that $%
\Delta _{j-\ell }(x)$ is a common ancestor of $\Delta _{j}(x)$ and its two
adjacent intervals. Since%
\begin{eqnarray*}
C_{n}\mu (\Delta _{j-\ell }(x)) &\leq &\mu (\Delta _{j}(x))\leq M_{j}(x) \\
&\leq &\mu (\Delta _{j}^{+}(x)\cup \Delta _{j}(x)\cup \Delta _{j}^{-}(x)) \\
&\leq &3\mu (\Delta _{j-\ell }(x))\text{ for all }l\leq L_{n}\text{,}
\end{eqnarray*}%
it will be sufficient to study the behaviour of the subsequences 
\begin{equation*}
\left\Vert T(\gamma _{0},\dots,\gamma _{J},\underbrace{\theta
_{1}^{-},\dots,\theta _{1}^{-}}_{k_{1}},\eta ^{o},\dots,\eta ^{\prime },%
\underbrace{\theta _{n+1}^{-},\dots,\theta _{n+1}^{-}}_{p_{n}})\right\Vert
\end{equation*}%
for $p_{n}\leq k_{n+1}$ and this we do in the same manner as in \cite{HHM}.
\end{proof}

It was shown in \cite[Thm. 5.7]{HHM} that the set of local dimensions at
essential points was a closed interval. Here we prove the same conclusion
for the set of local dimensions at truly essential points.

\begin{theorem}
\label{thm:convex} Suppose $\mu $ is a self-similar measure of finite type
satisfying the positive row property. Let $y,z$ be interior essential,
positive, periodic points. Then the set of local dimensions of $\mu $ at
truly essential points contains the closed interval with endpoints $\dim
_{loc}\mu (y)$ and $\dim _{loc}\mu (z)$.
\end{theorem}

\begin{proof}
Let $y$ and $z$ have truly essential, positive periods $\phi $ and $\theta $
respectively, with $T(\phi )=A$ and $T(\theta )=B$. Let $\eta _{1},\eta _{2}$
be truly essential, positive paths joining the last letter of $\phi $ to the
first letter of $\theta $ and vice versa. Given $0<t<1$, choose subsequences 
$m_{k},n_{k}\rightarrow \infty $ such that%
\begin{equation*}
\frac{L(\theta ^{-})m_{k}}{L(\theta ^{-})m_{k}+L(\phi ^{-})n_{k}}\rightarrow
t\text{.}
\end{equation*}%
Put%
\begin{equation*}
T(\psi _{k})=B^{m_{k}}T(\eta _{1})A^{n_{k}}T(\eta _{2})
\end{equation*}%
and consider a truly essential, positive, periodic point $x_{k}$ with period 
$\psi _{k}$. Using Lemma \ref{compare} we deduce that%
\begin{equation*}
sp(B^{m_{k}}T(\eta _{1})A^{n_{k}}T(\eta _{2}))\sim
sp(B)^{m_{k}}sp(A)^{n_{k}}.
\end{equation*}%
Coupled with Proposition \ref{periodic}, this implies%
\begin{eqnarray*}
\lim_{k}\dim _{loc}\mu (x_{k}) &=&\lim_{k}\frac{\log sp(B)^{m_{k}}+\log
sp(A)^{n_{k}}}{(L(\theta ^{-})m_{k}+L(\phi ^{-})n_{k})\log \varrho } \\
&=&t\dim _{loc}\mu (z)+(1-t)\dim _{loc}\mu (y).
\end{eqnarray*}%
Now appeal to the previous theorem.
\end{proof}

The three theorems combine to yield the following important corollary.

\begin{corollary}
\label{closedint} Let $\mu $ be a self-similar measure of finite type
satisfying the positive row property. Let $I=\inf \{\dim _{loc}\mu (x):x$
interior essential, positive, periodic$\}$ and $S=\sup \{\dim _{loc}\mu
(x):x $ interior essential, positive, periodic$\}.$ Then%
\begin{equation*}
\{\dim _{loc}\mu (x):x\text{ interior essential}\}=\{\dim _{loc}\mu (x):x%
\text{ truly essential}\}=[I,S]\text{.}
\end{equation*}%
A similar statement holds for the lower and upper local dimensions.
\end{corollary}

It is worth commenting here that this need not be the case for the set of
upper local dimensions of a maximal loop class (outside of the truly
essential class). An example is given in Section \ref{sec:MLC}.

\begin{remark}
In Example \ref{ex:new}, the local dimension of the boundary point $x$ with
symbolic representations $(1,4a,4a,4a,\dots )$ and $(1,3c,3d,3d,\dots )$ is 
\begin{equation*}
\frac{\left\vert \log (\max sp(T(4a,4a)),sp(T(3d,3d)))\right\vert }{\log 4}.
\end{equation*}%
Regardless of the choice of probabilities, this local dimension is always
contained within the interval that is the set of local dimensions of truly
essential points. It would be interesting to know if there were any examples
of self-similar measures of finite type and essential points $x$ where $\dim
_{loc}\mu (x)$ is not contained in the set of local dimensions of the truly
essential points.
\end{remark}

\subsection{Local dimension and the dimension of the support.}

In this section we show that, assuming the positive row property, the
essential class must contain a point $x$ such that $\dim _{loc}\mu (x)=\dim
_{H}(K)$.

\begin{lemma}
\label{tech}Let $\mu $ be a self-similar measure of finite type, with $%
s=\dim_{H}K$. Let $E$ denote the set of truly essential points and put 
\begin{align*}
G_{1}& =\{x\in E\mid \overline{\dim}_{loc}\mu (x)>s\}, \\
G_{2}& =\{x\in E\mid \underline{\dim}_{loc}\mu (x)<s\}.
\end{align*}%
Then $\mu (G_{1})=0=H^{s}(G_{2})$.
\end{lemma}

\begin{proof}
We recall that there are only countably many boundary essential points and
every non-atomic measure assigns mass zero to a countable set. Hence the
statement will be true if and only if it is true for $E$= the set of
interior essential points.

Let $x\in G_{1}$, say with $\overline{\dim }_{loc}\mu (x)=s(1+\varepsilon )$
for some $\varepsilon >0$. Then there will exist infinitely many $n$ such
that 
\begin{equation*}
\frac{\log \mu (\Delta _{n}(x))}{n\log \varrho }>s\left( 1+\frac{\varepsilon 
}{2}\right) .
\end{equation*}%
By Lemma \ref{lem:H(Delta)} we have 
\begin{equation*}
\lim_{n}\frac{\log H^{s}(\Delta _{n}(x))}{n\log \varrho }=s.
\end{equation*}%
This implies that there are infinitely many $n$ such that 
\begin{equation*}
\frac{\log \mu (\Delta _{n}(x))}{\log \varrho }>\left( 1+\frac{\varepsilon }{%
3}\right) \frac{\log H^{s}(\Delta _{n}(x))}{\log \varrho }
\end{equation*}%
and therefore 
\begin{equation*}
\mu (\Delta _{n}(x))\leq H^{s}(\Delta _{n}(x))H^{s}(\Delta
_{n}(x))^{\varepsilon /3}.
\end{equation*}%
Since $H^{s}(\Delta _{n})\rightarrow 0$ as $n\rightarrow \infty ,$ it
follows that for all $0<a<1$ there exists an $n$ such that 
\begin{equation*}
\mu (\Delta _{n}(x))\leq aH^{s}(\Delta _{n}(x)).
\end{equation*}%
In a similar way, if $x\in G_{2}$ and $b>1,$ then there exists an $n$ such
that 
\begin{equation*}
\mu (\Delta _{n}(x))\geq bH^{s}(\Delta _{n}(x)).
\end{equation*}

Define 
\begin{align*}
E_{1,n}^{a}& =\bigcup \{\Delta \in \mathcal{F}_{n}:\mu (\Delta )\leq
aH^{s}(\Delta )\},\ \ \ \ E_{1}^{a}=\bigcup\limits_{n}E_{1,n}^{a} \\
E_{2,n}^{b}& =\bigcup \{\Delta \in \mathcal{F}_{n}:\mu (\Delta )\geq
bH^{s}(\Delta )\},\ \ \ \ E_{2}^{b}=\bigcup\limits_{n}E_{2,n}^{b}
\end{align*}%
The comments above show that $G_{1}\subseteq E_{1}^{a}$ for all $0<a<1$ and $%
G_{2}\subseteq E_{2}^{b}$ for all $b>1$. Put 
\begin{equation*}
F_{1,1}^{a}=E_{1,1}^{a}\text{, }F_{1,n}^{a}=E_{1,n}^{a}\setminus \cup
_{k=1}^{n-1}F_{1,k}^{a}
\end{equation*}%
and similarly define $F_{2,n}^{b}$. Then $E_{1}^{a}$ is the disjoint union
of the sets $F_{1,n}^{a}$ and similarly for $E_{2}^{b}$. Further, we observe
that each set $F_{1,n}^{a}$ is a union of intervals, $\Delta ,$ with
disjoint interiors and the property that $\mu (\Delta )\leq aH^{s}(\Delta )$%
. Hence $\sigma $-additivity and the continuity of $H^{s}$ implies 
\begin{eqnarray*}
\mu (G_{1}) &\leq &\mu (E_{1}^{a})\leq \sum_{n}\mu (F_{1,n}^{a})\leq
a\sum_{n}H^{s}(F_{1,n}^{a}) \\
&=&aH^{s}(E_{1}^{a})\leq aH^{s}(K)\leq a.
\end{eqnarray*}%
As $0<a<1$ is arbitrary, we have that $\mu (G_{1})=0$.

Similarly 
\begin{equation*}
H^{s}(G_{2})\leq H^{s}(E_{2}^{b})\leq \frac{1}{b}\mu (E_{2}^{b})\leq \frac{1%
}{b}\mu (K)\leq \frac{1}{b}
\end{equation*}%
and as $b>1$ is arbitrary, we have that $H^{s}(G_{2})=0$.
\end{proof}

\begin{theorem}
\label{=Hdim}Let $\mu $ be a self-similar measure of finite type satisfying
the positive row property. Then there exists a truly essential element $x$
with $\dim _{loc}\mu (x)=\dim _{H}K$.
\end{theorem}

\begin{proof}
In fact, we will show a stronger result, that there exists an interior
essential point $x$ such that $\dim _{loc}\mu (x)=\dim _{H}K$.

Let $E$ be the set of interior essential points. According to Corollary \ref%
{closedint}, the set of local dimensions at the interior essential points is
an interval. So it suffices to show that the supremum of this interval is at
least $\dim _{H}K$ and the infimum is at most $\dim _{H}K$.

Assume, for a contradiction, that the infimum is strictly greater than $\dim
_{H}K$. This implies for all $x\in E$, 
\begin{equation*}
\dim _{H}K<\underline{\dim }_{loc}\mu (x)\leq \overline{\dim }_{loc}\mu (x)
\end{equation*}
and hence $E\subseteq G_{1}$. This fact, combined with Proposition \ref%
{SizeofTE} and Lemma \ref{tech}, gives $1=\mu (E)\leq \mu (G_{1})=0$, a
contradiction.

Similarly, if the supremum of the local dimensions of $E$ was strictly less
than $\dim _{H}K,$ then $E\subseteq G_{2}$, and hence $1=H^{s}(E)\leq
H^{s}(G_{2})=0$, a contradiction.
\end{proof}

It would be interesting to know if the set of such points has full $\mu $
measure. Notice that Lemma \ref{tech} implies that this is true if $\mu $ is
absolutely continuous with respect to $H^{s}$. Our next result gives
conditions under which this latter statement is true.

\begin{proposition}
\label{abscont}Suppose $\mu $ is a self-similar measure of finite type, with 
$\dim _{H}K=s$. Assume that the norm of any product of $n$ essential,
primitive transition matrices is bounded above by $C\varrho
^{sn(1-\varepsilon _{n})}$, where $\sup_{n}n\varepsilon _{n}<\infty $ and $%
C>0$ is a constant. Then $\mu $ is absolutely continuous with respect to $%
H^{s}$.
\end{proposition}

\begin{proof}
By \cite[p.35]{Ma}, $\mu \ll H^{s}$ if and only if $D(x)<\infty $ for $\mu $
almost all $x,$ where 
\begin{equation*}
D(x)=\liminf_{r\rightarrow 0}\frac{\mu (B(x,r))}{H^{s}(B(x,r))}.
\end{equation*}%
Appealing to Proposition \ref{SizeofTE}, we see it suffices to prove $%
D(x)<\infty $ for all interior essential points $x$. 
Standard arguments show it will be sufficient to prove 
\begin{equation*}
\liminf_{n\rightarrow \infty }\frac{\mu (\Delta _{n}(x)\cup \Delta
_{n}^{+}(x)\cup \Delta _{n}^{-}(x))}{H^{s}(\Delta _{n}(x)\cup \Delta
_{n}^{+}(x)\cup \Delta _{n}^{-}(x))}<\infty .
\end{equation*}

For $x$ an interior essential point, choose $J$ such that for all $n$
sufficiently large, $\Delta _{J+n}(x),\Delta _{J+n}^{+}(x)$ and $\Delta
_{J+n}^{-}(x)$ have a common essential ancestor at level $J$. 
Thus for $\Delta _{n+J}^{\prime }(x)$ denoting
any of $\Delta _{n+J}(x)$ or its two adjacent net intervals, 
\begin{eqnarray*}
\mu (\Delta _{n+J}^{\prime }(x)) &\sim &\left\Vert T(\gamma _{0},\dots
,\gamma _{J},\gamma _{J+1}^{\prime },\dots ,\gamma _{n+J}^{\prime
})\right\Vert \\
&\sim &\left\Vert T(\gamma _{J},\gamma _{J+1}^{\prime }\dots ,\gamma
_{n+J}^{\prime })\right\Vert \leq C\varrho ^{sn(1-\varepsilon _{n})}
\end{eqnarray*}%
for a constant $C$ not dependent on $n$. Here the last inequality comes from
the hypothesis of the proposition. Since Lemma \ref{lem:H(Delta)} implies $%
H^{s}(\Delta _{n})\sim \varrho ^{ns}$ for any $n$'th level net interval, 
\begin{equation*}
\liminf_{n\rightarrow \infty }\frac{\mu (\Delta _{n}(x)\cup \Delta
_{n}^{+}(x)\cup \Delta _{n}^{-}(x))}{H^{s}(\Delta _{n}(x)\cup \Delta
_{n}^{+}(x)\cup \Delta _{n}^{-}(x))}\leq \liminf_{n\rightarrow \infty }C%
\frac{\varrho ^{sn(1-\varepsilon _{n})}}{\varrho ^{ns}}<\infty
\end{equation*}%
as $\sup n\varepsilon _{n}<\infty $.
\end{proof}

\begin{remark}
We note that this proposition did not require the assumption of the positive
row property. Moreover, similar arguments show that $H^{s}|_{\mathrm{supp}%
\mu }$ is absolutely continuous with respect to $\mu $ if the norm of any
product of $n$ essential, primitive transition matrices is bounded below by $%
C\varrho ^{sn(1+\varepsilon _{n})}$, where $\sup_{n}n\varepsilon _{n}<\infty 
$ and $C>0$ is a constant.
\end{remark}

In the next example the self-similar measure is mutually absolutely
continuous to Lebesgue measure restricted to \textrm{supp}$\mu $ and the
local dimension is identical at all the truly essential points.

\begin{example}
\label{ex:abs cont} Consider the example $S_{j}(x)=x/4+b_{j}/12$ where $%
b_{j}\in \{0,1,2,7,8,9\}$ and associate to these the probabilities $%
p_{0}=p_{1}=p_{4}=p_{5}=1/8$, $p_{2}=p_{3}=1/4$. This measure does not have
full interval support, although the support is still of dimension one. To
see this, we observe that $K=[0,5/12]\cup \lbrack 7/12,1]$. There is one
reduced characteristic vector within the essential class. The four
transition matrices from this vector to itself are: 
\begin{equation*}
\left[ 
\begin{array}{ccc}
\frac{1}{8} & 0 & 0 \\ 
0 & 0 & \frac{1}{4} \\ 
\frac{1}{8} & \frac{1}{4} & 0%
\end{array}%
\right] ,\ \ \left[ 
\begin{array}{ccc}
\frac{1}{8} & \frac{1}{8} & 0 \\ 
0 & 0 & 0 \\ 
\frac{1}{8} & \frac{1}{8} & \frac{1}{4}%
\end{array}%
\right] ,\ \ \left[ 
\begin{array}{ccc}
\frac{1}{4} & \frac{1}{8} & \frac{1}{8} \\ 
0 & 0 & 0 \\ 
0 & \frac{1}{8} & \frac{1}{8}%
\end{array}%
\right] ,\ \ \left[ 
\begin{array}{ccc}
0 & \frac{1}{4} & \frac{1}{8} \\ 
\frac{1}{4} & 0 & 0 \\ 
0 & 0 & \frac{1}{8}%
\end{array}%
\right]
\end{equation*}%
We notice that all column sums of all of these matrices are exactly the same
at $1/4$. Hence the norm of any $n$-fold product of these matrices is
comparable to $4^{-n}$. This gives that the local dimension at all truly
essential points is $1$ and the measure $\mu $ is mutually absolutely
continuous with respect to Lebesgue measure on its support. It is worth
observing that this is true, despite this example not satisfying the
positive row property. Note that the points outside the essential class do
not necessarily have local dimension 1. For instance $\dim _{loc}\mu
(0)=\log 8/\log 4=3/2$.
\end{example}

Another illustration of this is seen in Example \ref{singleton} where this
phenomena occurs for a Cantor-like measure, when $H^{s}$ is the normalized
Lebesgue measure.

\section{Biased Bernoulli convolutions with simple Pisot contractions}

\label{sec:Pisot}

In this section we will assume $\mu $ is a Bernoulli convolution generated
by the IFS 
\begin{equation*}
\{S_{0}(x)=\varrho x,S_{1}(x)=\varrho x+(1-\varrho )\}
\end{equation*}%
and probabilities $p,1-p$, where $\varrho $ is the inverse of a simple Pisot
number (one whose minimal polynomial is of the form $x^{k}-x^{k-1}-\cdot
\cdot \cdot -x-1$) and $0<p<1$. The self-similar set is $[0,1]$, hence the
positive row property holds for all these Bernoulli convolutions.

Feng in \cite{F1} showed that if $p=1/2$, then $\mu $ has no isolated point
in its multifractal spectrum. In contrast, we will show here that if $p\neq
1/2$ there is always an isolated point, either $\dim _{loc}\mu (0)$ or $\dim
_{loc}\mu (1),$ depending on whether $p$ is less than or greater than $1/2$.

In \cite[Sect. 5]{F1}, Feng determined the characteristic vectors,
transition graph, and primitive transition matrices for the case $p=1/2$.
Using this information, it is not difficult to determine the primitive
transition matrices for the general case. In what follows, we use Feng's
notation to label the characteristic vectors as $a,b,d$, $c_{i}$, $\overline{%
c_{1}}$, $e_{j}$ $f_{j}$, $g,$ where $i=1,\dots,k$ and $j=1,\dots,k-1$. Here
all but $a,b,d$ are in the essential class.

\begin{lemma}
\label{lem:prim} The primitive transition matrices for the vectors in the
essential class are given by:

\begin{equation*}
\begin{array}{ll}
T(c_{j-1},c_{j})=%
\begin{bmatrix}
p & 0 \\ 
0 & 1-p%
\end{bmatrix}
\text{ for }2\leq j\leq k\text{,} & T(c_{k},g)=%
\begin{bmatrix}
p \\ 
1-p%
\end{bmatrix}
\text{,} \\ 
T(c_{k},c_{1})=%
\begin{bmatrix}
p & 0 \\ 
1-p & p%
\end{bmatrix}
\text{,} & T(c_{k},\overline{c_{1}})=%
\begin{bmatrix}
1-p & p \\ 
0 & 1-p%
\end{bmatrix}
\text{,} \\ 
T(g,f_{1})=T(f_{j},f_{j+1})=%
\begin{bmatrix}
p%
\end{bmatrix}
\text{,} & T(f_{j},c_{1})=%
\begin{bmatrix}
1-p & p%
\end{bmatrix}
\text{ for } j\leq k-2, \\ 
T(g,e_{1})=T(f_{j},e_{1})=%
\begin{bmatrix}
1-p%
\end{bmatrix}%
\text{,} & T(e_{j},f_{1})=%
\begin{bmatrix}
p%
\end{bmatrix}
\text{ for }j\leq k-1 \\ 
T(g,c_{1})=T(e_{j},c_{1})=%
\begin{bmatrix}
1-p & p%
\end{bmatrix}
\text{,} & T(e_{j},e_{j+1})=%
\begin{bmatrix}
1-p%
\end{bmatrix}
\text{ for } j\leq k-2.%
\end{array}%
\end{equation*}
\end{lemma}

\begin{proof}
We leave this as an exercise for the reader as it follows in a straight
forward manner from the information gathered in \cite{F1}. The main points
to observe are that if the $i$'th neighbour of a parent coincides with the $%
j $'th neighbour of a child, then $T_{ij}=p$, while if they differ by
(common) normalized distance $1-\varrho $, then $T_{ij}=1-p$. We also remind
the reader that for simple Pisot numbers, $\rho ^{-1},$ with minimal
polynomial of degree $k$, $1-\varrho =\varrho -\varrho ^{k}$.

We illustrate this with $T(c_{j-1},c_{j}).$ From \cite{F1} it can be seen
that the (normalized) neighbours of $c_{j}$ are $0$ and $1-\varrho ^{k-j+1},$
and $c_{j}$ is the only child of the parent $c_{j-1}$. If we renormalize so
they can be compared, we see that the two $0$ neighbours coincide and the
non-$0$ neighbours differ by $1-\varrho $. Thus $T$ is diagonal with the
entries being $p$ and $1-p$, respectively.
\end{proof}

\textbf{Notation:} Given a matrix $T$, denote by $\left\Vert T\right\Vert
_{\min }$ the pseudo-norm%
\begin{equation*}
\left\Vert T\right\Vert _{\min }=\min_{j}\sum_{i}\left\vert
T_{ij}\right\vert \text{ }
\end{equation*}%
where the sum is over all the rows of the matrix. That is, $\left\Vert T
\right\Vert_{\min}$ is the minimal column sum of $T$. 
Obviously, $\left\Vert T\right\Vert \geq \left\Vert T\right\Vert _{\min }$.
A useful property is that $\left\Vert T_{1}T_{2}\right\Vert _{\min }\geq
\left\Vert T_{1}\right\Vert _{\min }\left\Vert T_{2}\right\Vert _{\min }$.

\begin{lemma}
There exists an integer $N$ such that if $x\in (0,1)$ then%
\begin{equation*}
\lbrack x]=(\gamma _{1},\dots,\gamma _{M},\eta _{1},\eta _{2},\dots)
\end{equation*}
where $\gamma _{1},\dots,\gamma _{M}$ are characteristic vectors, $\eta _{j}$
are essential paths of length at most $N$ whose first letter, denoted $\eta
_{j,1}$, equals $c_{1}$, $\overline{c_{1}}$ or $f_{1},$ and 
\begin{equation*}
\left\Vert T(\eta _{j},\eta _{j+1},\eta _{j+2,1})\right\Vert _{\min }\geq
\min (p^{L-1}(1-p),(1-p)^{L-1}p)
\end{equation*}%
where $L=L(\eta _{j},\eta _{j+1})$.
\end{lemma}

\begin{proof}
One can see from the transition maps given in \cite[Sec 5.1]{F1} that the
symbolic representation for any $x\in (0,1)$ begins either as $%
[x]=(a,c_{1},\dots)$ or $[x]=(a,\ast ,\ast ,\dots,\ast ,y,\dots)$ where $%
\ast $ denotes (all) $b^{\prime }s$ or $d^{\prime }s$ and $y$ is either $%
c_{1},e_{1}$ or $f_{1}$. In the case when $y=e_{1}$ the path must continue
as $(e_{1},\dots,e_{j-1},z)$ where $z=c_{1}$ or $f_{1}$ and $j\leq k-2,$ or
as $(e_{1},\dots,e_{k-1},f_{1})$. Whichever is the case, one can see that
each essential $x$ must eventually admit either a (first) $c_{1}$ or $f_{1}$%
. This will be the first letter of $\eta _{1}$. Now define $\eta _{j}$ to
begin with the $j^{\text{th}}$ occurrence of either $c_{1}$ (or $\overline{%
c_{1}}$ in Feng's notation) or $f_{1}$. We need to check that with this
construction the $\eta _{j}$ are paths of bounded length (independent of $x$%
) and have the required property on the pseudo-norm of the transition
matrices.

First, suppose a path $\eta _{j}$ begins with $c_{1}$. Then it must continue
as $(c_{1},\dots,c_{k})$. If $c_{k}$ is followed by $c_{1}$ (or $\overline{%
c_{1}}$), then we stop and take $(c_{1},\dots,c_{k})$ as $\eta _{j}$ having
length $k$. Otherwise $c_{k}$ is followed by $g$ and if that is followed by $%
c_{1}$ or $f_{1}$ then $\eta _{j}=(c_{1},\dots,c_{k},g)$ has length $k+1$.
The only other possibility is that $g$ is followed by $e_{1}$, but in that
case, as we saw above, the path will continue as $(e_{2},\dots,,e_{j})$ with 
$j\leq k-1$, before continuing with either $c_{1}$ or $f_{1}$ (necessarily
with $f_{1}$ if $j=k-1$). Such a path $\eta _{j}$ has length at most $%
k+1+k-1=2k$.

To summarize, the paths $\eta _{j}$ that begin with $c_{1},$ together with
the first letter of $\eta _{j+1},$ are of the form $(c_{1},\dots,c_{k})$
with next letter either $c_{1}$ or $\overline{c_{1}}$, $(c_{1},%
\dots,c_{k},g) $ with next letter either $c_{1}$ or $f_{1}$, or $%
(c_{1},\dots,c_{k},g,e_{1},\dots,e_{j})$ with $j\leq k-1$ and next letter
either $c_{1}$ or $f_{1}$ (necessarily $f_{1}$ if $j=k-1)$.

The arguments are similar for the paths that begin with $f_{1}$, with these
paths having length at most $2k-2.$

Now we verify the claimed pseudo-norm property. For this we apply the
previous lemma to analyze the product of the appropriate transition
matrices. Of course, any primitive transition matrix has pseudo-norm at
least $\min (p,1-p)$.

For paths $\eta _{J}$ that begin with $c_{1}$, we will see that even 
\begin{equation*}
\left\Vert T(\eta _{J},\eta _{J+1,1})\right\Vert _{\min }\geq \min \left(
p^{L-1}(1-p),(1-p)^{L-1}p\right) \text{ for }L=L(\eta _{J})
\end{equation*}%
and this will certainly imply the claim. To prove this, we consider the
different paths individually.

\textbf{Case 1}: $(\eta _{J},\eta _{J+1},_{1})=(c_{1},\dots,c_{k},y)$ with $%
y_{1}=\eta _{J+1,1}=c_{1}$ or $\overline{c_{1}}$: If $y=c_{1}$, then 
\begin{equation*}
T(\eta _{J},\eta _{J+1,1})=T(c_{1},\dots,c_{k})T(c_{k},c_{1})=\left[ 
\begin{array}{cc}
p^{k} & 0 \\ 
(1-p)^{k} & (1-p)^{k-1}p%
\end{array}%
\right]
\end{equation*}%
and hence has pseudo-norm with the required lower bound. The argument when
the first letter of $\eta _{J+1}=$ $\overline{c_{1}}$ is similar.

\textbf{Case 2:} $(\eta _{J},\eta _{J+1},_{1})=(c_{1},\dots,c_{k},g,y)$ with 
$y=c_{1}$ or $f_{1}$: If $y=c_{1}$, then an easy calculation shows 
\begin{equation*}
T(\eta _{J},\eta _{J+1,1})=T(c_{1},\dots,c_{k})T(c_{k},g)T(g,c_{1})=\left[ 
\begin{array}{cc}
p^{k}(1-p) & p^{k+1} \\ 
(1-p)^{k+1} & p(1-p)^{k}%
\end{array}%
\right] .
\end{equation*}%
If $y=\,f_{1}$, then 
\begin{equation*}
T(\eta _{J},\eta _{J+1,1})=\left[ 
\begin{array}{c}
p^{k+1} \\ 
(1-p)^{k}p%
\end{array}%
\right] .
\end{equation*}

The cases $(\eta _{J},\eta
_{J+1},_{1})=(c_{1},\dots,c_{k},g,e_{1},\dots,e_{j},y) $ for $y=c_{1}$ or $%
f_{1}, $ or $(\eta _{J},\eta
_{J+1},_{1})=(f_{1},\dots,f_{j},e_{1},\dots,e_{i},f_{1})$ are similar.

The only case where we must consider two consecutive paths, $\eta _{J}\eta
_{J+1},$ are when $\eta _{J}=(f_{1},\dots,f_{j})$ and either the next letter
is $c_{1}$ or the path continues as $(f_{1},\dots,f_{j},e_{1},\dots,e_{i})$
with $i,j\geq 1$ and the next letter is $c_{1}$. But in that case the next
path, $\eta _{J+1},$ is one of the paths beginning with $c_{1}$ discussed
above and we already know that then 
\begin{equation*}
\left\Vert T(\eta _{J+1},\eta _{J+2,1})\right\Vert _{\min }\geq \min
(p^{L-1}(1-p),(1-p)^{L-1}p)\text{ for }L=L(\eta _{J+1}).
\end{equation*}%
Combining this bound with the fact that $\left\Vert T(\eta _{J},\eta
_{J+1,1})\right\Vert _{\min }\geq \min (p^{L},(1-p)^{L})$ for $L=L(\eta
_{J}) $ completes the proof.
\end{proof}

\begin{theorem}
\label{prop:Simple Pisot} Suppose $\mu $ is a Bernoulli convolution with
contraction factor $\varrho $ the inverse of a simple Pisot number, and with
probabilities $p\neq 1-p$. Then there is an isolated point in the set of
local dimensions of $\mu $ at either $0$ or $1$, depending on which of $p$
or $1-p$ is smaller.
\end{theorem}

\begin{proof}
Without loss of generality assume $p<1/2$. Standard arguments show that $%
\dim_{loc}\mu (0)=\log p/\log $ $\varrho $.

Consider any $x\in (0,1)$. As the set of local dimensions of boundary
essential points is contained in the set of local dimensions of interior
essential points, we can assume without loss of generality that $x$ is an
interior essential point. Write $[x]=(\gamma _{1},\dots,\gamma _{M},\eta
_{1},\eta _{2},\dots)$ with the notation as in the previous lemma. We have
the formula 
\begin{equation*}
\dim_{loc}\mu (x)=\lim_{J}\frac{\log \left\Vert T(\eta _{1},\eta
_{2},\dots,\eta _{2J},\eta _{2J+1,1})\right\Vert }{\sum_{i=1}^{2J}L(\eta
_{i})\log \varrho },
\end{equation*}%
should the local dimension of $\mu $ at $x$ exist$.$

Set $L_{i}=L(\eta _{2i-1},\eta _{2i})$. Then 
\begin{eqnarray*}
\left\Vert T(\eta _{1},\eta _{2},\dots ,\eta _{2J},\eta
_{2J+1,1})\right\Vert &\geq &\prod_{i=1}^{J}\left\Vert T(\eta _{2i-1},\eta
_{2i},\eta _{2i+1,1})\right\Vert _{\min } \\
&\geq &p^{\sum_{i}(L_{i}-1)}(1-p)^{J}.
\end{eqnarray*}%
Hence%
\begin{equation*}
\log \left\Vert T(\eta _{1},\eta _{2},\dots ,\eta _{2J},\eta
_{2J+1,1})\right\Vert =\left( \sum_{i}L_{i}-J\right) \log p+J\log (1-p),
\end{equation*}%
so that 
\begin{equation*}
\frac{\log \left\Vert T(\eta _{1},\eta _{2},\dots ,\eta _{2J},\eta
_{2J+1,1})\right\Vert }{\sum_{i=1}^{J}L_{i}}\geq \log p+\frac{J(\log
(1-p)-\log p)}{\sum_{i=1}^{J}L_{i}}\text{.}
\end{equation*}%
But $L_{i}\leq 2N$ (where $N$ is as in the lemma), hence for any $J$, 
\begin{equation*}
\frac{\log \left\Vert T(\eta _{1},\eta _{2},\dots ,\eta _{2J},\eta
_{2J+1,1})\right\Vert }{\sum_{i=1}^{2J}L(\eta _{i})\log \varrho }\leq \frac{%
\log p}{\log \varrho }+\frac{\log (1-p)-\log p}{2N\log \varrho }<\frac{\log p%
}{\log \varrho }=\dim _{loc}\mu (0).
\end{equation*}%
and therefore $\dim _{loc}\mu (x)$ is bounded away from $\dim _{loc}\mu (0)$.
\end{proof}

\section{Cantor-like measures of Finite type}

\label{sec:Cantor}

The focus of this section will be the Cantor-like self-similar sets and
measures generated by the IFS 
\begin{equation}
\left\{ S_{j}(x)=\frac{1}{d}x+\frac{j}{md}(d-1):j=0,\dots ,m\right\}
\label{Cantor}
\end{equation}%
for integers $d\geq 2$ and probabilities $p_{j}>0$, $j=0,\dots ,m.$ The
self-similar set is the $m$-fold sum of the Cantor set with contraction
factor $1/d,$ rescaled to $[0,1],$ and is the full interval when $m\geq d-1$%
. We will assume this to be the case for otherwise the IFS satisfies the
open set condition and is well understood. This class of measures includes,
for example, the $m$-fold convolution of the uniform Cantor measure
associated with the Cantor set generated by $S_{0}(x)=\frac{1}{d}x$, $%
S_{1}(x)=\frac{1}{d}x+\frac{d-1}{d}$. As $K=[0,1]$ when $m\geq d-1$, we see
that all of these examples satisfy the positive row property.

These measures were studied by different methods in \cite{BHM} and \cite{Sh}
where it was shown, for example, that $\{\dim_{loc}\mu (x):x\in (0,1)\}$ was
a closed interval. In \cite[Sect. 7]{HHM} it was shown that the essential
class for any of these Cantor-like measures is $(0,1),$ hence all $x\in
(0,1) $ are truly essential. Consequently, the fact that $\{\dim _{loc}\mu
(x):x\in (0,1)\}$ is a closed interval can also be deduced from our
Corollary \ref{closedint}.

In this section we will establish more refined information about the local
dimensions of these measures. In particular, we give a new proof of the fact
that $\dim_{loc}\mu (0)$ (or $\dim_{loc}\mu (1)$) is an isolated point if $%
p_{0}$ (resp., $p_{m}$) is the minimal probability, as was shown by other
methods in \cite{BHM} and \cite{Sh}. We give an example to show that there
need not be an isolated point if this is not the case, as well as examples
of Cantor-like measures whose set of local dimensions consists of
(precisely) two points.

For this detailed analysis it is helpful to completely determine the finite
type structure of these measures. There are two cases to consider, $m\equiv
0 $ mod$(d-1)$ and $m\neq 0$ mod$(d-1).$

\begin{proposition}
\label{thm:equiv d-1} Assume $\mu $ is the self-similar Cantor-like measure
of finite type generated by the IFS (\ref{Cantor}), with $m=k(d-1)$ for
integer $k$.

\begin{enumerate}
\item The essential class has one reduced characteristic vector, $E$, with
normalized length $1/k$ and neighbour set $(j/k:j=0,\dots ,k-1)$. The
reduced characteristic vector $E$ has $d$ children, identical to itself,
labelled as $E^{(i)}$, $i=1,\dots ,d$. \label{part:1}

\item There are $m-k+2$ net intervals at level one with reduced
characteristic vector $E$. These are the intervals $[\frac{k-1}{kd},\frac{k}{%
kd}],\dots,[1-\frac{k}{kd},1-\frac{k-1}{kd}]$. \label{part:2}

\item The primitive transition matrix $T(E,E^{(i)})$ is given by the
following formula: For $x,y=0,..,k-1$, 
\begin{equation*}
(T(E,E^{(i)}))_{x,y}=\left\{ 
\begin{array}{ll}
p_{dx-y+i-1} & \text{if}\ 0\leq dx-y+i-1\leq m \\ 
0 & \text{otherwise}%
\end{array}%
\right. .
\end{equation*}%
\label{part:3}
\end{enumerate}
\end{proposition}

\begin{example}
Consider the IFS as in (\ref{Cantor}) with $d=4$ and $m=9$, $k=3$. The
essential class consists of the one reduced characteristic vector $%
(1/3,(0,1/3,2/3))$. There are four transition matrices from $E$ to $E$. They
are 
\begin{equation*}
\left[ 
\begin{array}{ccc}
p_{0} & 0 & 0 \\ 
p_{4} & p_{3} & p_{2} \\ 
p_{8} & p_{7} & p_{6}%
\end{array}%
\right] \text{, }\left[ 
\begin{array}{ccc}
p_{1} & p_{0} & 0 \\ 
p_{5} & p_{4} & p_{3} \\ 
p_{9} & p_{8} & p_{7}%
\end{array}%
\right] \text{, }\left[ 
\begin{array}{ccc}
p_{2} & p_{1} & p_{0} \\ 
p_{6} & p_{5} & p_{4} \\ 
0 & p_{9} & p_{8}%
\end{array}%
\right] \text{, }\left[ 
\begin{array}{ccc}
p_{3} & p_{2} & p_{1} \\ 
p_{7} & p_{6} & p_{5} \\ 
0 & 0 & p_{9}%
\end{array}%
\right] .
\end{equation*}
\end{example}

\begin{proof}[Proof of Proposition \protect\ref{thm:equiv d-1}.]
As noted in the proof of Proposition 7.1 of \cite{HHM}, 
\begin{equation*}
\{S_{\sigma }(0):\sigma \in \mathcal{A}^{n}\}=\left\{ \frac{(d-1)j}{md^{n}}%
:0\leq j\leq (d^{n}-1)k\right\}
\end{equation*}%
and 
\begin{equation*}
\{S_{\sigma }(1):\sigma \in \mathcal{A}^{n}\}=\left\{ \frac{(d-1)(j+k)}{%
md^{n}}:0\leq j\leq (d^{n}-1)k\right\} .
\end{equation*}%
First, consider the level $n$ net intervals that lie in $[1/d^{n},1-1/d^{n}]$%
. These have the form 
\begin{equation*}
\Delta ^{(j)}=\left[ \frac{(d-1)j}{d^{n}m},\frac{(d-1)(j+1)}{d^{n}m}\right] =%
\text{ }\left[ \frac{j}{d^{n}k},\frac{j+1}{d^{n}k}\right]
\end{equation*}%
for $j=k,\dots,k(d^{n}-1)-1$. They have normalized length $(d-1)/m=1/k$ and
normalized neighbours as claimed in \eqref{part:1} of the statement of the
Proposition. These net intervals have $d$ children, 
\begin{equation*}
\left[ \frac{(d-1)(dj+i)}{d^{n+1}m},\frac{(d-1)(dj+i+1)}{d^{n+1}m}\right] 
\text{ for }i=0,\dots,d-1\text{,}
\end{equation*}%
all of the same type again.

At level $1$, the net intervals have the form $[\frac{j-1}{dk},\frac{j}{dk}]$%
. If $j<k-1$, then there are only $j$ neighbours because $j-1-i<0$ if $i\geq
j$. If $j>m$ there are $<k$ neighbours because $(j-1)/dk$ is not an iterate
of $0$. All other net intervals are type $E$. This proves \eqref{part:2}.

Now consider the $x$ neighbour of $E^{(j)}$ at level $n$, for $0\leq x\leq
k-1,$ namely 
\begin{equation*}
S_{\sigma _{x}}(0)=\frac{(d-1)(j-x)}{d^{n}m},\text{ }
\end{equation*}%
and the $y$ neighbour of its $i$'th child, $E^{(i)}$, for $0\leq y\leq k-1$,%
\begin{equation*}
S_{\sigma _{y}}(0)=\frac{(d-1)(dj+i-1-y)}{d^{n+1}m}.
\end{equation*}%
For any $0\leq w\leq m$ it follows that 
\begin{equation*}
S_{\sigma _{x}w}(0)=\frac{(d-1)(d(j-x)+w)}{d^{n+1}m}.
\end{equation*}%
Hence, whenever $0\leq dx-y+i-1=w\leq m$ we have $S_{\sigma
_{x}w}(0)=S_{\sigma _{y}}(0)$, and this proves \eqref{part:3}.
\end{proof}

\begin{example}
\label{1/2}Suppose $m=k(d-1)$ is even. Then $1/2=S_{\sigma }(0)$ for some $%
\sigma \in \mathcal{A}$ and therefore $1/2$ is a left endpoint of a net
interval of level one, and hence is a boundary essential point. Thereafter, $%
1/2$ is the left endpoint of the left-most child of the parent net interval
and thus $1/2$ has symbolic representation $(E^{(m/2)},E^{(1)},E^{(1)},\dots
)$. Similarly, $1/2$ is also the right-most endpoint of the right-most child
of the net interval immediately to the left of this net interval.
Consequently, $1/2$ also has symbolic representation $%
(E^{(m/2)-1},E^{(d)},E^{(d)},\dots )$.

When $k=2$ ($m=2(d-1)$), for example, then $T(E^{(1)},E^{(1)})=\left[ 
\begin{array}{cc}
p_{0} & 0 \\ 
p_{d} & p_{d-1}%
\end{array}%
\right] $ and $T(E^{(d)},E^{(d)})=\left[ 
\begin{array}{cc}
p_{d-1} & p_{d-2} \\ 
0 & p_{m}%
\end{array}%
\right] $, so that we have 
\begin{equation*}
\dim_{loc}\mu (1/2)=\left\vert \log (\max (p_{0},p_{d-1},p_{m}))\right\vert
/\log d.
\end{equation*}
\end{example}

\begin{proposition}
\label{thm:nequiv d-1} Assume $\mu $ is the self-similar Cantor-like measure
of finite type generated by the IFS (\ref{Cantor}), with $m=k(d-1)+r$, $%
1\leq r\leq d-2$.

\begin{enumerate}
\item The essential class consists of two reduced characteristic vectors, $E$
with normalized length $r/m$ and neighbour set $(j(d-1)/m$ $:j=0,\dots,k)$,
and $F$ with normalized length $(d-1-r)/m$ and neighbour set $%
((r+j(d-1))/m:j=0,\dots,k-1)$.\label{part:1a}

\item At level one the essential net intervals are alternately $E$ and $F$,
beginning with the interval $[\frac{1}{d}-\frac{r}{md},\frac{1}{d}]$ of type 
$E$ and ending with $[1-\frac{1}{d},1-(\frac{1}{d}-\frac{r}{md})]$ also type 
$E$. There are $m-k+1$ net intervals with characteristic vector $E$ and $m-k$
with characteristic vector $F$.

\item Type $E$ has $2r+1$ children labelled (from left to right) $E^{(1)}$, $%
F^{(2)}$, \dots, $E^{(2r+1)}$.

Type $F$ has $2(d-r)-1$ children labelled $F^{(1)}$, $E^{(2)}$, \dots,$%
F^{2(d-r)-1}$.\label{part:2a}

\item The non-zero entries of the primitive transition matrices are as
follows:
\end{enumerate}

\begin{itemize}
\item For $i=0,\dots,r$ and $0\leq x,y\leq k$, $%
(T(E,E^{(2i+1)}))_{xy}=p_{dx-y+i}$ if $\ 0\leq dx-y+i\leq m;$

\item For $i=1,\dots,r$ and $0\leq x\leq k,0\leq y\leq k-1$, $%
(T(E,F^{(2i)}))_{xy}=p_{dx-y+i-1}$ if $0\leq dx-y+i-1\leq m;$

\item For $i=0,\dots,d-r-1$ and $0\leq x,y\leq k-1$, $%
(T(F,F^{(2i+1)}))_{xy}=p_{dx+r-y+i}$ if $0\leq dx+r-y+i\leq m;$

\item For $i=1,\dots,d-r-1$ and $0\leq x\leq k-1,0\leq y\leq k$, $%
(T(F,E^{(2i)}))_{xy}=p_{dx+r-y+i}$ if $0\leq dx+r-y+i\leq m.$\label{part:3a}
\end{itemize}
\end{proposition}

\begin{proof}
The proof is similar to the previous case, but with two characteristic
vectors arising because the iterates of $0$ and $1$ do not coincide. Indeed%
\begin{equation*}
\{S_{\sigma }(1):\sigma \in \mathcal{A}^{n}\}=\left\{ \frac{(d-1)(j+k)+r}{%
md^{n}}:0\leq j\leq (d^{n}-1)k\right\} .
\end{equation*}%
The net intervals whose left endpoint is an iterate of $0$ give one
characteristic vector and those whose left endpoint is an iterate of $1$ is
the second. We leave the details for the reader.
\end{proof}

\begin{example}
\label{1/2r}Suppose $k=1$, $m=d-1+r$ where $1\leq r\leq d-2$ and $m$ is
even. There are an odd number of net intervals at level one and by symmetry $%
1/2$ lies at the centre of the middle interval. This is a net interval of
type $F$ since $2(m-1)+1\equiv 3\mod 4$. At all other levels there are an
odd number of net intervals, so again $1/2$ lies at the centre of the middle
one and again this is a type $F,$ namely $F^{(2i+1)}$ where $i=(d-r-1)/2$,
since $2(d-r)-1\equiv 1$mod $4$. As $T(F^{(2i+1)},F^{(2i+1)})=[p_{m/2}]$, we
have $\dim_{loc}\mu (1/2)=\left\vert \log p_{m/2}\right\vert /\log d$.
\end{example}

In the proof of the next result we will use the pseudo norm $\left\Vert
T\right\Vert _{\min }$, defined in the previous section, and also the norm 
\begin{equation*}
\left\Vert T\right\Vert _{\max }\text{ =}\max_{j}\sum_{i}\left\vert
T_{ij}\right\vert
\end{equation*}%
where the sum is over all the rows of the matrix. That is, $\left\Vert T
\right\Vert_{\max}$ is the maximal column sum of $T$. 
For matrices with non-negative values it is easy to see that 
\begin{equation*}
\left\Vert T_{1}T_{2}\right\Vert _{\min }\geq \left\Vert T_{1}\right\Vert
_{\min }\left\Vert T_{2}\right\Vert _{\min }, \ \ \left\Vert
T_{1}T_{2}\right\Vert _{\max }\leq \left\Vert T_{1}\right\Vert _{\max
}\left\Vert T_{2}\right\Vert _{\max }
\end{equation*}%
and 
\begin{equation*}
\left\Vert T\right\Vert _{\min }\leq \left\Vert T\right\Vert \leq
C\left\Vert T\right\Vert _{\max }
\end{equation*}%
where $C$ is the number of columns of $T$.

\begin{proposition}
\label{thm:Cantor} Let $P_{i}=\sum_{i\equiv j\mod d}p_{j}$, $P_{\max }=\max
(P_{i})$, and $P_{\min }=\min (P_{i})$. For any $x\in (0,1)$, we have 
\begin{equation*}
\frac{\left\vert \log P_{\max }\right\vert }{\log d}\leq \underline{\dim}%
_{loc} \mu (x)\leq \overline{\dim}_{loc}\mu (x)\leq \frac{\left\vert \log
P_{\min }\right\vert }{\log d}.
\end{equation*}
\end{proposition}

\begin{proof}
From the formulas given in Proposition \ref{thm:equiv d-1} and \ref%
{thm:nequiv d-1} one can see that the column sums of an essential, primitive
transition matrix $T$ are of the form $P_{i}$. Hence if $T$ is a product of $%
m$ essential, primitive transition matrices, then 
\begin{equation*}
P_{\min }^{m}\leq \left\Vert T\right\Vert \leq CP_{\max }^{m}
\end{equation*}%
where $C$ is a bound for the number of columns of a primitive transition
matrix.

Since any $x\in (0,1)$ is truly essential and the set of local dimensions of
boundary essential points is contained in the set of local dimensions of
interior essential points (Cor. \ref{closedint}), we can assume without loss
of generality that $x$ is an interior essential point. Hence there exists a $%
k$ so that $\Delta _{k}(x)$ is an essential net interval and a common
ancestor for $\Delta _{n}^{-}(x)$, $\Delta _{n}(x)$ and $\Delta _{n}^{+}(x)$
for all $n\geq k$. Consequently, $\mu (\Delta _{n}^{-}(x)),\mu (\Delta
_{n}(x))$ and $\mu (\Delta _{n}^{+}(x))$ can all be approximated by the
norms of products of $n-k$ primitive transition matrices within the
essential class. From this the result follows.
\end{proof}

\begin{corollary}
\begin{enumerate}
\item If $p_{0}$ $<P_{\min }$, then $\dim_{loc}\mu (0)$ is an isolated
point. \label{co:1}

\item If $m\geq d$ and $p_{0}$ $<p_{j}$ for $j\neq 0,m$, then $\dim_{loc}\mu
(0)$ is an isolated point. \label{co:2}

Similar statements can be made for $p_{m}$ and $\dim_{loc}\mu (1)$.
\end{enumerate}
\end{corollary}

\begin{proof}
We have that \eqref{co:1} is immediate since $\dim_{loc}\mu (0)=|\log
p_{0}|/\log d$.

For \eqref{co:2}, one can easily check from these formulas that $p_{0}$ (and 
$p_{m}$) are never the only non-zero entries in a column when $m\geq d$.
Hence the hypothesis of \eqref{co:1} is satisfied.
\end{proof}

\begin{remark}
We remark that it is possible for \eqref{co:1} to be satisfied without $%
p_{0} $ being minimal. For instance, if $m\geq 2d$, then every column admits
at least two non-zero entries and hence it would suffice to have $%
p_{0}<2p_{j}$ for all $j $ in order for $\dim_{loc}\mu (0)$ to be an
isolated point.
\end{remark}

On the other hand, it is also possible for such a measure to have no
isolated points. Here is an example.

\begin{example}
Consider the IFS $\{S_{j}(x)=x/3+j/6$ $:j=0,\dots ,4\}$ and probabilities $%
p_{0}=p_{4}=1/3$, $p_{1}=p_{2}=p_{3}=1/9$. The essential class is composed
of one reduced characteristic vector, with three transition matrices from
this vector to itself. The transition matrices are 
\begin{equation*}
\begin{bmatrix}
1/3 & 0 \\ 
1/9 & 1/9%
\end{bmatrix}%
\text{, }%
\begin{bmatrix}
1/9 & 1/3 \\ 
1/3 & 1/9%
\end{bmatrix}%
\text{, }%
\begin{bmatrix}
1/9 & 1/9 \\ 
0 & 1/3%
\end{bmatrix}%
.
\end{equation*}%
One can check that the second matrix has $4/9$ as an eigenvalue. Further,
all matrices have maximal column sum equal to $4/9$. This gives an exact
lower bound for the set of local dimensions. One can compute that the local
dimension of the essential class, $(0,1)$, contains the interval $I=[\frac{%
\log (9/4)}{\log 3},1.24]\approx \lbrack 0.738,1.24]$ and is contained in $[%
\frac{\log (9/4)}{\log 3},2.00]$. We can establish the upper bounds by
explicitly finding a point of local dimension 1.24 in the first case, and by
using the $||\cdot ||_{max}$ norm for the second case. The local dimension
of the self-similar measure at the two end points of the support is $1$ and $%
1\in I$.
\end{example}

\begin{corollary}
If $P_{\max }=P_{\min },$ then $\{\dim_{loc}\mu (x):x\in (0,1)\}$ $=\{1\}$.
\end{corollary}

\begin{proof}
This follows from the observation that $d\cdot P_{\min }\leq \sum P_{i}\leq
d\cdot P_{\max }$ and $\sum P_{i}=1$.
\end{proof}

Here is a family of examples of this phenomena, generalizing \cite[Ex. 6.1]%
{HHM}.

\begin{example}
\label{singleton}Suppose $\mu $ is the self-similar measure associated to
the IFS (\ref{Cantor}) with $m+1\equiv 0$ mod $d$ and $p_{j}=1/(m+1)$ for
all $j=0,\dots,m\geq d$. Then 
\begin{equation*}
\dim_{loc}\mu (x)=1\text{ for all }x\in (0,1)
\end{equation*}%
and%
\begin{equation*}
\dim_{loc}\mu (0)=\dim_{loc}\mu (1)=\frac{\log (m+1)}{\log d}>1,
\end{equation*}%
so the set of local dimensions is a doubleton.
\end{example}

\begin{proof}
The assumption that $m+1\equiv 0$ mod $d$ ensures that each column of each
essential primitive transition matrix $T$ has exactly $k$ non-zero entries,
where $m+1=kd$. Consequently, $P_{i}=k/(m+1)=1/d$ for each $i$ and the
result follows from the previous corollary.
\end{proof}

\begin{remark}
These measures are also an example of the phenomena addressed in Proposition %
\ref{abscont}. The proof above shows there exists a constant $C$ such that $%
d^{-n}\leq \left\Vert T\right\Vert \leq Cd^{-n}$for all $n$-fold products of
primitive transition matrices. As $\dim _{H}$ supp$\mu =1$, the proposition
implies $\mu $ restricted to the truly essential class is absolutely
continuous with respect to Lebesgue measure on $[0,1]$.
\end{remark}

\begin{corollary}
Suppose $\{\mu _{n}\}$ is a sequence of Cantor-like measures, all with
contraction factor $1/d$. Let $P_{\max }^{(n)}$ and $P_{\min }^{(n)}$ be the
maximal and minimal column sums associated with $\mu _{n}$. If $P_{\max
}^{(n)}$ $-P_{\min }^{(n)}\rightarrow 0$, then the set of local dimensions
at any $x\in (0,1)$ tends to $1$.
\end{corollary}

\begin{proof}
Similar reasoning to the proof of the previous corollary shows that $P_{\min
}^{(n)},P_{\max }^{(n)}\rightarrow \frac{1}{d}$.
\end{proof}

\begin{example}
Let $\mu $ be the self-similar measure associated to the IFS (\ref{Cantor})
and let $\mu ^{k}$ be the $k$-fold convolution of $\mu $, normalized to $%
[0,1]$. Then 
\begin{equation*}
\dim_{loc}\mu ^{k}(x)\rightarrow 1\text{ for all }x\in (0,1)\text{ and }%
\dim_{loc}\mu ^{k}(x)\rightarrow \infty \text{ for }x=0,1.
\end{equation*}

To see this, let $Q(x)=p_{0}+p_{1}x+\dots +p_{n}x^{n}$. The measure $\mu
^{k} $ is also a Cantor-like measure with contraction factor $1/d$. With the
contractions ordered in the natural way, the probability of the $j^{th}$
term, denoted $p_{j}^{(k)},$ is equal to the coefficient of $x^{j}$ in $%
Q(x)^{k}$ and 
\begin{equation*}
P_{i}^{(k)}=\sum_{j\equiv i\text{ }mod\text{ }d}p_{j}^{(k)}=\frac{1}{d}%
\sum_{j=1}^{d}Q(\zeta _{d}^{j})^{k}\zeta _{d}^{-ji}
\end{equation*}%
where $\zeta _{d}$ is a primitive $d^{th}$ root of unity. It is easy to see
that $Q(1)\ =1$ and $|Q(\zeta _{d}^{j})|<1$ for $j\neq d$. Hence we see that 
$P_{i}^{(k)}\rightarrow 1/d$ as $k\rightarrow \infty $ for all $i$. This in
turn implies that $P_{\min }^{(k)}-P_{\max }^{(k)}\rightarrow 0$ and hence $%
\dim_{loc}\mu ^{k}(x)\rightarrow 1$ for each $x\in (0,1)$.

In contrast, $\dim_{loc}\mu ^{k}(0)=\lim \left\vert \log
p_{0}^{k}\right\vert /\log d\rightarrow \infty $ and similarly for $\dim
_{loc}\mu ^{k}(1).$
\end{example}

\begin{example}
Suppose $\nu $ is the uniform Cantor measure associated with the IFS $%
\,\{S_{0}(x)=x/d$, $S_{1}(x)=x/d+(d-1)/d\}$. Then $\nu ^{m}$ is the measure
of finite type generated by the IFS (\ref{Cantor}) and probabilities $%
p_{j}=2^{-m}\binom{m}{j}$. Information was given about the minimum and
maximum local dimensions (other than at $0,m$) in \cite[Thm. 6.1]{BHM} for $%
m\leq 2d-1$.

We can extend the maximum local dimension result to $m<3(d-1)$ when $m-d$ is
odd, as follows. First, note that the column sums of essential primitive
transition matrices have the form%
\begin{equation*}
P_{j}=2^{-m}\sum_{k=-\infty }^{\infty }\binom{m}{j+kd}
\end{equation*}%
and reasoning as in \cite[Lem.6.2]{BHM} shows that these are minimized at $j=%
\left[ \frac{m-d}{2}\right] $. We can assume $m=2(d-1)+r$ for $1\leq r\leq
d-2$. Consider the periodic element $x_{0}$ with period $\theta
=(F^{(2i+1)},F^{(2i+1)})$ for $i=(d-r-1)/2=(m-d+1)/2-r$. Then 
\begin{equation*}
T(F^{(2i+1)},F^{(2i+1)})=\left[ 
\begin{array}{cc}
p_{\frac{m-d+1}{2}} & p_{\frac{m-d-1}{2}} \\ 
p_{\frac{m+d+1}{2}} & p_{\frac{m+d-1}{2}}%
\end{array}%
\right] .
\end{equation*}%
The two column sums are equal and minimal among all column sums of essential
primitive transition matrices. Thus $\left\Vert T\right\Vert \sim \left\Vert
T\right\Vert _{\min }$ and further, this is a lower bound on the norm of any
essential primitive transition matrix. Hence $\dim _{loc}\mu (x_{0})$ is
maximal over all $x\in (0,1)$.

Since the column sums are maximized when $j=[m/2]$, we deduce from Examples %
\ref{1/2} and \ref{1/2r} that $\dim _{loc}\nu ^{m}(1/2)=\left\vert \log
p_{m/2}\right\vert /\log d$ is minimal when $m<2(d-1)$ is even, as was also
seen in \cite{BHM}.
\end{example}

\section{Maximal loop classes outside the essential class}

\label{sec:MLC}

In \cite{HHM}, it is shown that if $\mu $ is a self-similar measure of
finite type, with full support and regular probabilities, then the set of
upper (or lower) local dimensions at points in any maximal loop class is an
interval. In this section we show that this is not true for finite type
measures satisfying only the positive row property. The example we use is a
self-similar measure that would be Cantor-like, in the sense of the previous
section, if there we had allowed some probabilities to be zero.

The measure $\mu $ will arise from the maps $S_{i}(x)=x/4+d_{i}/12$ with $%
d_{i}=i$ for $i=0,1,\dots,5$, $d_{6}=8$ and $d_{7}=9$, and probabilities $%
p_{0}=1/2$, $p_{i}=1/14$ for $i=1,\dots,7$. The reduced transition diagram
has 7 reduced characteristic vectors. The reduced characteristic vectors are:

\begin{itemize}
\item Reduced characteristic vector 1: $(1, (0))$

\item Reduced characteristic vector 2: $(1/3, (0))$

\item Reduced characteristic vector 3: $(1/3, (0, 1/3))$

\item Reduced characteristic vector 4: $(1/3, (0, 1/3, 2/3))$

\item Reduced characteristic vector 5: $(1/3, (1/3, 2/3))$

\item Reduced characteristic vector 6: $(1/3, (2/3))$

\item Reduced characteristic vector 7: $(2/3, (0, 1/3))$
\end{itemize}

The maps are:

\begin{itemize}
\item RCV $1 \to [2, 3, 4, 4, 4, 4, 5, 6, 2, 7, 6]$

\item RCV $2 \to [2, 3, 4, 4]$

\item RCV $3 \to [4, 4, 4, 4]$

\item RCV $4 \to [4, 4, 4, 4]$

\item RCV $5 \to [4, 4, 5, 6]$

\item RCV $6 \to [2, 7, 6]$

\item RCV $7 \to [4, 4, 4, 4, 4, 4, 5, 6]$
\end{itemize}

\begin{figure}[tbp]
\includegraphics[scale=0.5]{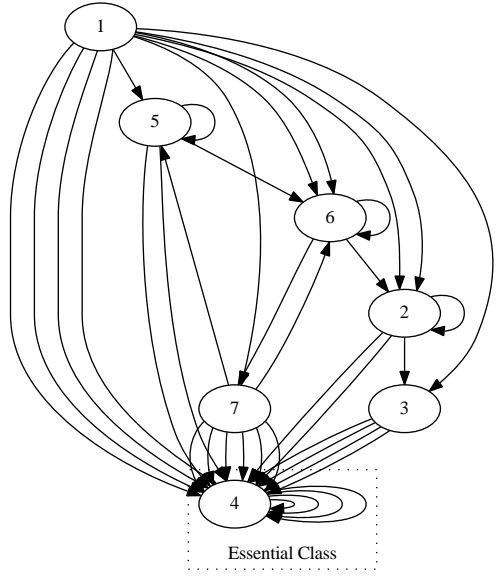}
\caption{Transition diagram for example in Section \protect\ref{sec:MLC}}
\label{fig:Pic3}
\end{figure}

We refer to Figure \ref{fig:Pic3} for the transition diagram. 

\begin{figure}[tbp]
\includegraphics[scale=0.45]{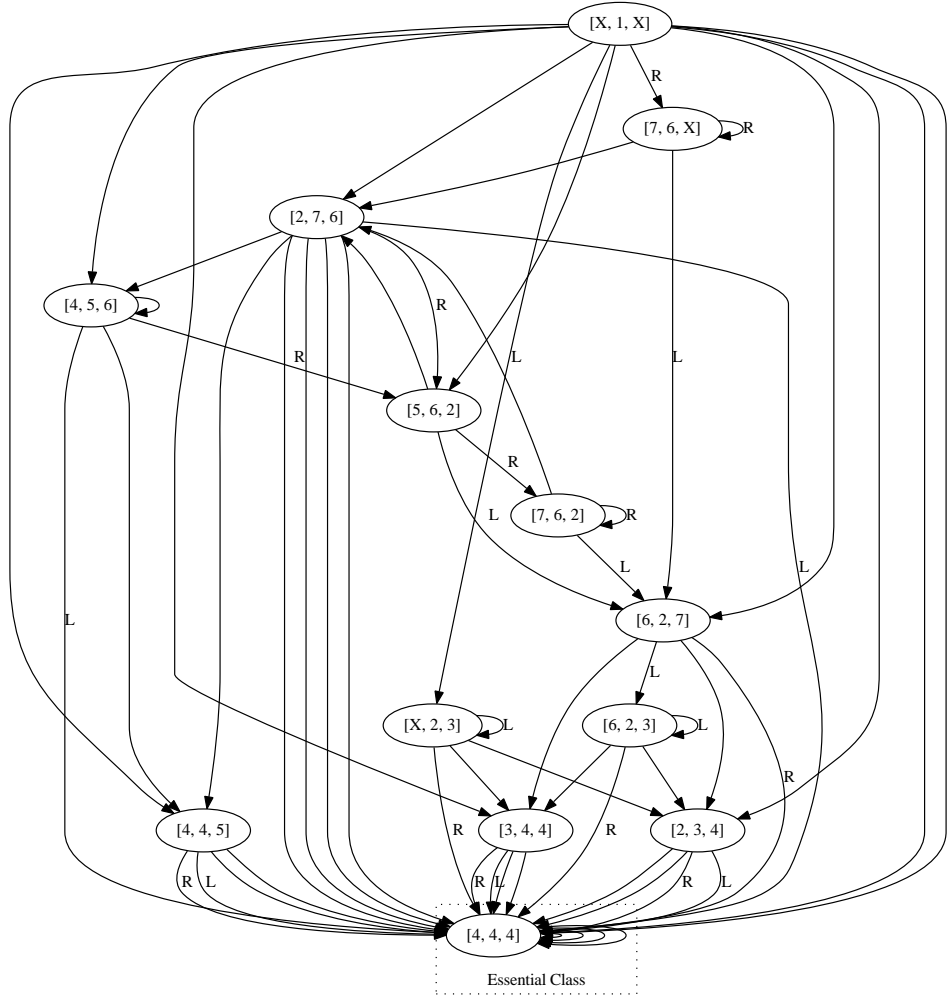}
\caption{Triple transition diagram for example in Section \protect\ref%
{sec:MLC}}
\label{fig:Pic4}
\end{figure}

As the probabilities are not regular, the reduced transition diagram does
not contain all of the necessary information to compute the local dimension
at a point, since to calculate $\dim_{loc} \mu (x)$ we need to know about $%
\Delta _{n}^{-}(x),$ $\Delta _{n}^{+}(x),$ in addition to $\Delta _{n}(x)$.
To keep track of this information, we introduce the triple transition
diagram. Each triple consists of a net interval and its adjacent net
intervals. If there is no adjacent net interval, then we represent this with
an $X$. The triple transition diagram also displays the transitions from
each triple to their triple children and denotes which transitions are right
or left-most descendents. See Figure \ref{fig:Pic4} for the triple
transition diagram.

We define, in the obvious way, the triple loop classes, triple maximal loop
classes and the triple essential class. In this example, the set of points
that are in the triple essential class, $[4,4,4]$, is the same as the set of
truly essential points. To see this, observe that if $x$ is an interior
essential point, then there exists an integer $n$ such that $x$ is in the
interior of the net interval $\Delta _{n}(x)$ whose reduced characteristic
vector is of type 4. As $x$ is not equal to the end point of $\Delta _{n}(x)$%
, there will exist some $k$ such that $\Delta _{n+k}(x)$ and its two
neighbours will all have reduced characteristic vector of type $4$. Hence $%
[\Delta _{n+k}^{-}(x),\Delta _{n+k}(x),\Delta _{n+k}^{+}(x)]=[4,4,4]$. If,
instead, $x$ is a boundary essential point, then there exists an $n$ such
that two adjacent $\Delta _{n}(x)$ and $\Delta _{n}^{\prime }(x)$ are the
reduced characteristic vector of type $4$. In this case, regardless of which
net interval containing $x$ we use, we see that $\Delta _{n+1}(x)$ and its
two adjacent net intervals will be the reduced characteristic vector of type 
$4$. Hence $[\Delta _{n+1}^{-}(x),\Delta _{n+1}(x),\Delta
_{n+1}^{+}(x)]=[4,4,4]$. The other inclusion is clear.

\begin{figure}[tbp]
\includegraphics[scale=0.5]{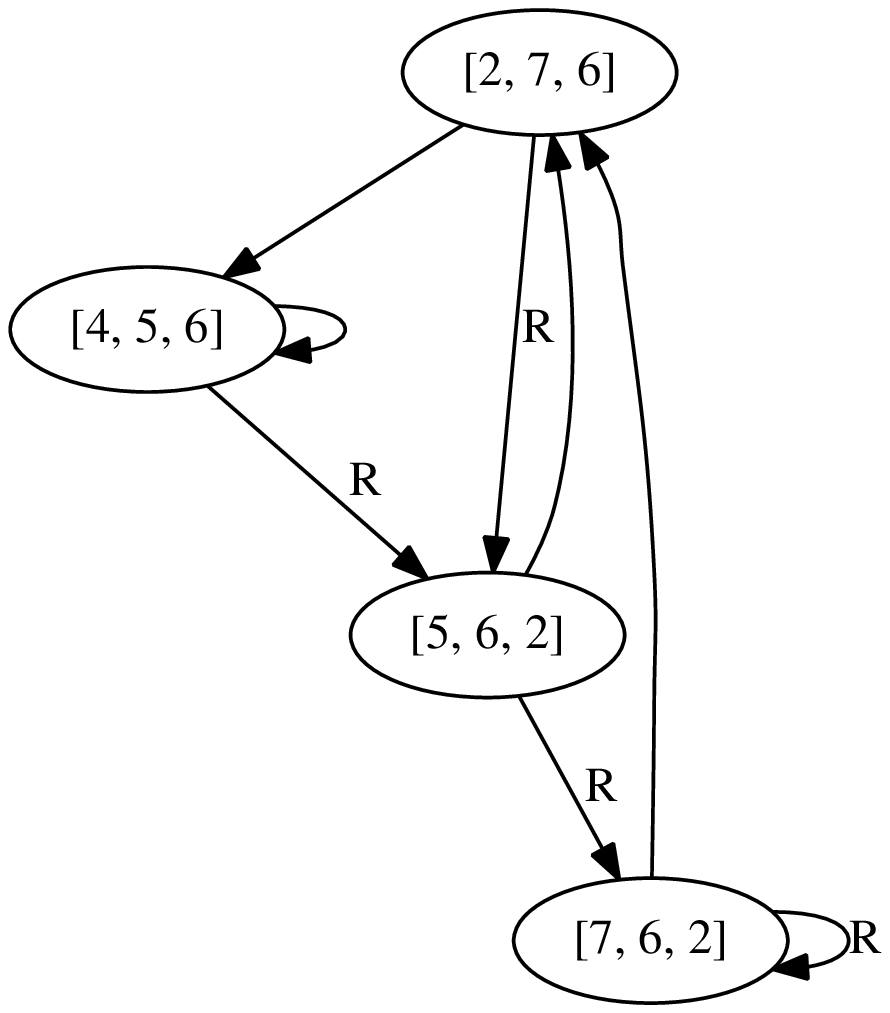}
\caption{Triple transition diagram for maximal loop class for example in
Section \protect\ref{sec:MLC}}
\label{fig:Pic5}
\end{figure}

From the triple transition diagram, we can see that there are four triple
maximal loop classes, in addition to the triple essential class. Three of
these are singletons, $[7,6,X]$, $[X,2,3]$ and $[6,2,3]$. It is very easy to
compute the local dimensions of these points. The final maximal loop class
is formed by the four triples $[2,7,6],[4,5,6],[5,6,2],[7,6,2]$ and is of
positive type. See Figure \ref{fig:Pic5} for the triple transition diagram
of this triple maximal loop class. We have indicated on this diagram which
of these transitions are right or left-most descendents.

We will determine the local dimension of points in this (non-singleton)
triple maximal loop. It is important to note that this triple loop class
admits no left-most descendents.

First, assume that the symbolic representation of a point $x$ in the loop
class does not contain arbitrarily long, right-most paths, say these lengths
are bounded by $K$. This implies that $\Delta _{n}(x)$ is in the interior of 
$\Delta _{n-K}(x)$, hence $\Delta _{n}(x),\Delta _{n}^{+}(x)$ and $\Delta
_{n}^{-}(x)$ are all comparable to $\Delta _{n-K}(x)$ for all $n$. Thus we
may ignore the $\Delta _{n}^{+}(x)$ and $\Delta _{n}^{-}(x)$ and this allows
us to use the techniques from \cite{HHM} without modification. (We will not
be able to ignore $\Delta _{n}^{+}(x)$ and $\Delta _{n}^{-}(x)$ later, when
we allow arbitrarily long right-most paths.)

In this case, the relevant transition matrices are: 
\begin{align*}
T(5,5)& =T(7,5)=\left[ 
\begin{array}{cc}
1/14 & 1/14 \\ 
1/14 & 1/14%
\end{array}%
\right] & T(5,6)& =T(7,6)=\left[ 
\begin{array}{c}
1/14 \\ 
1/14%
\end{array}%
\right] \\
T(6,7)& =\left[ 
\begin{array}{cc}
1/14 & 1/14%
\end{array}%
\right] & T(6,6)& =\left[ 
\begin{array}{c}
1/14%
\end{array}%
\right]
\end{align*}

For these matrices, the minimal column sum is $1/14$, and the maximal sum is 
$1/7$. These numbers are also the eigenvalues of $T(6,6)$ and $T(5,5)$
respectively. As we are only concerned with $\Delta_n(x)$, and do not need
to worry about $\Delta_n^+(x)$ or $\Delta_n^{-}(x)$, we see that the
standard convexity argument can be used to show that the set of local
dimensions is an interval. Consequently, such points produce the interval 
\begin{equation*}
\left[ \frac{\log 7}{\log 4},\frac{\log 14}{\log 4}\right] \approx \lbrack
1.403677461,1.903677461]
\end{equation*}
as the set of local dimensions.

To consider the the case when $x$ contains arbitrarily long right-most
paths, we now need to consider $\Delta^+_n(x)$ and $\Delta_n^-(x)$. We will
need to know about the additional transition matrices

\begin{align*}
T(2,2) & = \left[ 
\begin{array}{c}
1/2%
\end{array}
\right] & T(6,2) & = \left[ 
\begin{array}{c}
1/14%
\end{array}
\right] \\
& & 
\end{align*}

First, consider an $x$ whose tail consists of the right-most branch of the
triples $[7,6,2],[7,6,2],[7,6,2],\dots $. We observe in this case that $%
\Delta _{n}^{-}(x)$ is comparable to $\Delta _{n}(x)$ as they share the
common ancestor $\Delta _{n-1}(x)$, so that $\mu (\Delta _{n}(x))$ and $\mu
(\Delta _{n}^{-}(x))$ are comparable to $\left\Vert T(6,6)^{n}\right\Vert
=1/14^{n}$. We further see that $\Delta _{n}^{+}(x)$ is \textbf{not}
comparable to $\Delta _{n}(x)$, as it does not share a common ancestor a
bounded number of generation back. In fact, the symbolic representation of $%
\Delta _{n}^{+}(x)$ has tail $(2,2,2,\dots ,2)$ and hence $\mu (\Delta
_{n}^{+}(x))$ is comparable to $\left\Vert T(2,2)^{n}\right\Vert =1/2^{n}$.
This gives us that the local dimension at $x$ is 
\begin{align*}
\dim _{loc}\mu (x)& =\lim_{n}\frac{\log M_{n}(x)}{n\log 1/4} \\
& =\lim_{n}\frac{\log \left( \mu (\Delta _{n}^{-}(x))+\mu (\Delta
_{n}(x))+\mu (\Delta _{n}^{+}(x)\right) }{n\log 1/4} \\
& =\lim_{n}\frac{\log \left( (1/14)^{n}+(1/14)^{n}+(1/2)^{n}\right) }{n\log
1/4} \\
& =1/2
\end{align*}

Next, consider the case where $x$ has arbitrarily long right-most paths, but
not infinitely-long, right-most paths from $[7,6,2]\rightarrow \lbrack
7,6,2] $. We claim that in this case the upper local dimension must be
greater than $\log 7/\log 4\sim 1.403677461$. To see this, we note that for
all $n$ where $\Delta _{n}(x)$ is not a right-most child of $\Delta
_{n-1}(x) $ (which happens infinitely often) the value of $\mu
(M_{n}(x))\sim \mu (\Delta _{n}(x))$, as $\mu (\Delta _{n}(x))\sim \mu
(\Delta _{n}^{+}(x))\sim \mu (\Delta _{n}^{-}(x))$. As on this subsequence
we have that the lim sup must be greater than $\log 7/\log 4,$ it follows
that the set of upper local dimensions is not an interval.

This is in contrast to the lower local dimension, where we can achieve any
value $z$ in the interval $\left[ \frac{1}{2},\frac{\log 14}{\log 4}\right] $%
. We will prove this by constructing an $x$ in this maximal loop class such
that \underline{dim}$_{loc}\mu (x)=z$.

Let $A=T((7,6,2),(2,7,6))\cdot T((2,7,6),(5,6,2))\cdot T((5,6,2),(7,6,2))$
be a triple of the transition matrices for the path through $%
(7,6,2)\rightarrow (2,7,6)\rightarrow _{R}(5,6,2)\rightarrow _{R}(7,6,2)$.
We note here that these transition matrices may work on the middle or the
right-most matrix of the previous transition, depending upon the nature of
the transition. Consider the path with transition matrices 
\begin{equation*}
T_{k}:=T((7,6,2),(7,6,2))^{n_{1}}\cdot A\cdot
T((7,6,2),(7,6,2))^{n_{2}}\cdot A\dots A\cdot T((7,6,2),(7,6,2))^{n_{k}}.
\end{equation*}%
We let $x$ be the point in $K$ with symbolic path $\lim_{k}T_{k}$. Let $L_{k}
$ be the length of $T_{k}$, that is, $L_{k}=n_{1}+2+n_{2}+2+\dots +2+n_{k}$.
We see that the three matrices associated with $T_{1}$ are 
\begin{equation*}
\left( 
\begin{bmatrix}
14^{-n_{1}} & 14^{-n_{1}}%
\end{bmatrix}%
,%
\begin{bmatrix}
14^{-n_{1}}%
\end{bmatrix}%
,%
\begin{bmatrix}
2^{-n_{1}}%
\end{bmatrix}%
\right) =\left( 
\begin{bmatrix}
14^{-L_{1}} & 14^{-L_{1}}%
\end{bmatrix}%
,%
\begin{bmatrix}
14^{-L_{1}}%
\end{bmatrix}%
,%
\begin{bmatrix}
2^{-L_{1}}%
\end{bmatrix}%
\right) .
\end{equation*}%
The three matrices associated with $T_{2}$ are 
\begin{align*}
& \left( 2%
\begin{bmatrix}
14^{-(n_{1}+2\ +n_{2})} & 14^{-(n_{1}+2+\ n_{2})}%
\end{bmatrix}%
,2%
\begin{bmatrix}
14^{-(n_{1}+2\ +n_{2})}%
\end{bmatrix}%
,%
\begin{bmatrix}
2^{-(n_{2}+1)}14^{-(n_{1}+1)}%
\end{bmatrix}%
\right) \\
& =\left( 2%
\begin{bmatrix}
14^{-L_{2}} & 14^{-L_{2}}%
\end{bmatrix}%
,2%
\begin{bmatrix}
14^{-L_{2}}%
\end{bmatrix}%
,%
\begin{bmatrix}
2^{-(n_{2}+1)}14^{-(L_{1}+1)}%
\end{bmatrix}%
\right) .
\end{align*}%
In general, for $k\geq 2$, we have that the three matrices associated to $%
T_{k}$ are 
\begin{equation*}
\left( 2^{k-1}%
\begin{bmatrix}
14^{-L_{k}} & 14^{-L_{k}}%
\end{bmatrix}%
,2^{k-1}%
\begin{bmatrix}
14^{-L_{k}}%
\end{bmatrix}%
,2^{k-2}%
\begin{bmatrix}
2^{-(n_{k}+1)}14^{-(L_{k-1}+1)}%
\end{bmatrix}%
\right)
\end{equation*}%
So, on the subsequence associated to $L_{k}$ we see that $M_{L_{k}}(x)$ is
approximately $2^{k-2}2^{-(n_{k}+1)}14^{-(L_{k-1}+1)}$. Choosing the $n_{k}$
such that 
\begin{equation*}
z=\lim_{k\rightarrow \infty }\frac{\log
(2^{k-2}2^{-(n_{k}+1)}14^{-(L_{k-1}+1)})}{\log (4^{-L_{k}})}
\end{equation*}%
gives that the local dimension, computing along this subsequence, is equal
to $z$. For example taking $n_{k}\approx \frac{\log 14-z\log 4}{(2z-1)\log 2}%
L_{k-1}$ will suffice. Note: so long as $z\in (\frac{1}{2},\frac{\log 14}{%
\log 4})$ we see that this is always a positive constant times $L_{k-1}$.

It is straightforward to see that this subsequence of lower local dimension
estimates gives a lower bound for the sequence, which proves the desired
result. To see this just note that if we consider a path for $x$ of length $%
N\in (L_{k-1},L_{k})$, then $M_{L_{k}}(x))^{1/L_{k}}<M_{N}(x))^{1/N}$.

Thus the set of lower local dimensions at points in the loop class is the
interval $\left[ \frac{1}{2},\frac{\log 14}{\log 4}\right]$. This is in
contrast to the set of upper local dimensions at points in the loop class,
which is the union of the interval $\left[ \frac{\log 7}{\log 4},\frac{\log
14}{\log 4}\right] $ together with the singleton $1/2$.

\section{When finite type IFS have Pisot contractions}

\label{sec:misc}

In this section we explore the connection between finite type and Pisot
contraction factors. This was motivated by Feng's observation in \cite{F4}
showing that the IFS $\{S_{j}(x)=\varrho x$ $+j(1-\varrho )/m:j=0,\dots,m\}$
satisfies the finite type condition if and only if $\varrho ^{-1}$ is Pisot.

In Example \ref{Ex:pos}, the IFS is of finite type, does not satisfy the
open set condition, but the contraction factor is not necessarily the
inverse of a Pisot number. This example also illustrates that we can have a
measure of finite type whose support is not the full interval $[0,1]$, yet
every row of each primitive transition matrix admits a non-zero entry. In
addition, it has the interesting property that every element of the
self-similar set is truly essential.

\begin{example}
\label{Ex:pos} Pick any positive number $\varepsilon <1/8$. Let $0<\varrho
<1 $ be a root of $\varepsilon -2x^{2}+4x-1$ and consider the self-similar
set $K$ generated by the contractions $S_{i}(x)=\varrho x+d_{i},$ with $%
d_{0}=0$, $d_{1}=-\varrho ^{2}+\varrho $, $d_{2}=\varepsilon -\varrho
^{2}+2\varrho $, and $d_{3}=\varepsilon -2\varrho ^{2}+3\varrho $. Consider
the associated probability measure with uniform probabilities, $p_{i}=1/4$
for $i=0,\dots ,3 $. There are 5 reduced characteristic vectors: $(1,(0))$, $%
(1-\varrho ,(0))$, $(\varrho ,(0,1-\varrho ))$, $(1-\varrho ,(\varrho ))$,
and $(1-2\varrho ,(\varrho ))$. Figure \ref{fig:nonzerorow} shows the
transition diagram. The essential class are all the characteristic vectors
except $1$ and there are no loop classes outside of the essential class.
Hence $K$ is the truly essential set. As this satisfies the positive row
property, the set of local dimensions is a closed interval.

We list below the transition matrices that are not equal to $[1/4]$. 
\begin{figure}[tbp]
\includegraphics[scale=0.7]{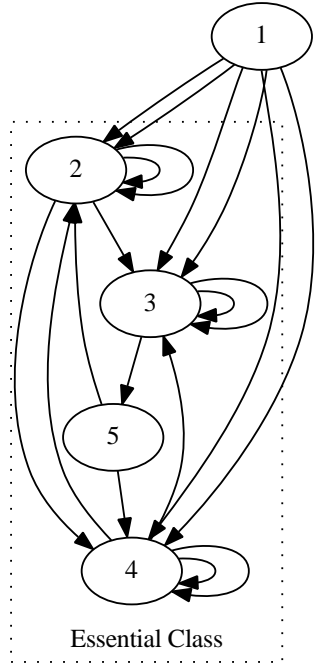}
\caption{Transition diagram for Example \protect\ref{Ex:pos}}
\label{fig:nonzerorow}
\end{figure}

\begin{equation*}
T(1,3)=T(2,3)=T(4,3)=\left[ 
\begin{array}{cc}
1/4 & 1/4%
\end{array}%
\right] \text{, }T(3,5)=\left[ 
\begin{array}{c}
1/4 \\ 
1/4%
\end{array}%
\right]
\end{equation*}%
\begin{equation*}
T(3,3)=\left[ 
\begin{array}{cc}
1/4 & 0 \\ 
1/4 & 1/4%
\end{array}%
\right] ,T(3,3)=\left[ 
\begin{array}{cc}
1/4 & 1/4 \\ 
0 & 1/4%
\end{array}%
\right] .
\end{equation*}%
Using techniques similar to \cite{HHM} one can show that the minimal local
dimension is 
\begin{equation*}
\frac{\log \left( sp%
\begin{bmatrix}
1/4 & 1/4 \\ 
0 & 1/4%
\end{bmatrix}%
\begin{bmatrix}
1/4 & 0 \\ 
1/4 & 1/4%
\end{bmatrix}%
\right) }{2\log \varrho }=\frac{\log \frac{3+\sqrt{5}}{32}}{2\log \varrho }
\end{equation*}%
and the maximal local dimension is%
\begin{equation*}
\frac{\log \left( sp%
\begin{bmatrix}
1/4 & 1/4 \\ 
0 & 1/4%
\end{bmatrix}%
\right) }{\log \varrho }=\frac{\log 1/4}{\log \varrho }.
\end{equation*}%
The details are left to the reader.

The incidence matrix of the essential class is 
\begin{equation*}
I=\left[ 
\begin{array}{cccc}
2 & 1 & 1 & 0 \\ 
0 & 2 & 0 & 1 \\ 
1 & 1 & 2 & 0 \\ 
1 & 0 & 1 & 0%
\end{array}%
\right] .
\end{equation*}%
Its spectral radius is $2+\sqrt{2}$, thus the formula from Proposition \ref%
{dimH} gives that $\dim_{H}K=\log (2+\sqrt{2})/\left\vert \log \varrho
\right\vert $.
\end{example}

In this example the overlap was `perfect', that is, all overlaps were of the
form $\varrho ^{n}K$ for some integer $n,$ ($\varrho $ the contraction
factor, $K$ the self-similar set). But $K\neq \lbrack 0,1]$. In our final
proposition we show that if the self-similar set is a full interval and the
overlaps are perfect, in this sense, then $\varrho $ is Pisot.

\begin{proposition}
Suppose $[0,1]$ is the self-similar set associated with contractions $S_{j},$
each with contraction factor $\varrho $. Assume that, for each $j$, the
length of the interval $S_{j}([0,1])\cap S_{j+1}([0,1]))$ is either equal to 
$\varrho ^{k_{j}}$ for some integer $k_{j}$ or has length equal to $0$. Then 
$\varrho $ is Pisot.
\end{proposition}

\begin{proof}
Assume that we have $n$ contractions. As the self-similar set is $[0,1]$, we
have that%
\begin{equation*}
n\varrho -\sum_{i=1}^{n-1}\varrho ^{k_{i}}=1.
\end{equation*}%
Let $N=\max(k_{j})$ and $q=\varrho ^{-1}.$ Let $f(z)=z^{N}-nz^{N-1}$ and $%
g(z)=\sum_{j}z^{N-k_{j}}$. Then $(f+g)(q)=0$. Clearly $f(n)=0$ and all other
zeros of $f$ are inside the unit disc (namely, at $0$). Further, on the unit
disc $\left\vert f(z)\right\vert \geq n-1\geq \left\vert g(z)\right\vert $.
By Rouche's theorem, $f+g$ has $n-1$ zeros in the closure of the unit disk
and therefore its other root, $q$, is a Pisot number.
\end{proof}

\begin{remark}
It would be interesting to fully understand the connection between finite
type and a Pisot contraction factor. Note that if $\dim _{H}K=1$, then as $%
\dim _{H}K=\log (sp(I))/\left\vert \log \varrho \right\vert $, and the
incidence matrix $I$ is integer valued, it follows that $\varrho ^{-1}$ is
an algebraic integer.
\end{remark}

\end{document}